\documentclass[12pt]{article}
\usepackage{amssymb}

\begin{document}
\newcommand{\lap}{\mbox{$\bigtriangleup$}}
\newcommand{\grad}{\mbox{$\bigtriangledown$}}
\newcommand{\be}{\begin{equation}}
\newcommand{\ee}{\end{equation}}
\newcommand{\kernel}{\frac{1}{ |x-y|^{n-\alpha} } }
\newcommand{\ind}{ (n+\alpha)/(n-\alpha) }
\newcommand{\indi}{ \frac{n+\alpha}{n-\alpha} }
\newcommand{\bound}{ (1+|x|)^{\alpha-n} }
\newcommand{\bounds}{\frac{1}{ (1+|x|)^{n-\alpha} } }
\newcommand{\solution}{ (\frac{t}{t^{2}+|x-x_{0}|^{2} } )^{ (n-\alpha)/2 } }
\newcommand{\sol}{ u_{t,x_{0}}(x) }
\newcommand{\half}{ \frac{1}{2} }
\newcommand{\hsp}{ \Sigma_{\lambda} }
\newcommand{\hspn}{ \Sigma_{\lambda}^{-} }
\newcommand{\hspp}{ \Sigma_{\lambda}^{+} }
\newcommand{\dom}{D(\lambda, \delta, R_{1})}
\newtheorem{mthm}{Theorem}
\newtheorem{mcor}{Corollary}
\newtheorem{mpro}{Proposition}
\newtheorem{mfig}{figure}
\newtheorem{mlem}{Lemma}
\newtheorem{mdef}{Definition}
\newtheorem{mrem}{Remark}
\newtheorem{mpic}{Picture}
\newtheorem{rem}{Remark}[section]
\newcommand{\ra}{{\mbox{$\rightarrow$}}}
\newtheorem{thm}{Theorem}[section]
\newtheorem{pro}{Proposition}[section]
\newtheorem{lem}{Lemma}[section]
\newtheorem{defi}{Definition}[section]
\newtheorem{cor}{Corollary}[section]

\title{\bf Semilinear equations involving the fractional Laplacian on domains}

\author{ Wenxiong Chen \thanks{Partially supported by the Simons Foundation Collaboration Grant for Mathematicians 245486.}
\hskip 0.5cm  Yanqin Fang \thanks{Partially supported by Young Teachers Program of Hunan University 602001003}
\hskip 0.5cm  Ray Yang \thanks{Partially supported by NSF Grant DMS-1103786}}
\date{\today}
\maketitle

\begin{abstract}
Let $\Omega$ be either a unit ball or a half space. Consider the following Dirichlet problem involving the fractional Laplacian for $0<\alpha<2$:
\begin{equation}
\left\{\begin{array}{ll}
(- \lap)^{\alpha/2} u = f(u), & x \in \Omega, \\
u \equiv 0 , & x \not{\in} \Omega .
\end{array}
\right.
\label{01}
\end{equation}
 
 Instead of using the conventional extension method of Cafarelli and Silvestre \cite{CaS}, we employ a new and direct approach by studying the equivalent integral equation
 \begin{equation}
 u(x) = \int_{\Omega} G(x, y) f(u(y)) d y ,
 \label{02}
 \end{equation}
 where $G(x,y)$ is the Green's function associated with the fractional Laplacian in the domain $\Omega$.
 Applying the {\em method of moving planes in integral forms}, we establish radial symmetry for
positive solutions in the unit ball and obtain Liouville type theorems--non-existence of positive solutions in the half space.
The regularity of solutions are also investigated.

\end{abstract}

\bigskip


\section{Introduction}

The fractional Laplacian in $\mathbb{R}^n$ is a nonlocal integral operator, taking the form
\begin{equation}
(-\Delta)^{\alpha/2} u(x) = C_{n,\alpha} P.V. \int_{\mathbb{R}^n} \frac{u(x)-u(z)}{|x-z|^{n+\alpha}} dz
\label{Ad7}
\end{equation}
 where $\alpha$ is any real number between $0$ and $2$ and P.V. means in the Cauchy principal value sense.

In recent years, there has been a great deal of interest in using the fractional Laplacian to model diverse physical phenomena, such as anomalous diffusion and quasi-geostrophic
flows, turbulence and water waves, molecular dynamics, and relativistic quantum mechanics
of stars ( see [BoG] [CaV] [Co] [TZ] and the
references therein). It also has various applications in probability
and finance [A] [Be] [CT]. In particular, the fractional Laplacian
can be understood as the infinitesimal generator of a stable L\'{e}vy process [Be].
 We refer the readers to Di Nezza, Palatucci, and Valdinoci's survey paper \cite{NPV} for a detailed exposition of the function spaces involved in the analysis of the operator and a long list of relevant references.

Let $\Omega$ be a domain in $\mathbb{R}^n$. In this paper, we analyze the behavior of solutions to the Dirichlet problem for semilinear equations
\begin{equation}
\left\{\begin{array}{ll}
(-\Delta)^{\alpha/2} u(x) = f(u(x)) , & x \in \Omega , \\
u(x) \equiv 0 , & x \not{\in} \Omega .
\end{array}
\right.
\label{00}
\end{equation}
We study symmetry, monotonicity, regularity, and non-existence of positive solutions.

There are several distinctly different ways to define the fractional Laplacian in a domain $\Omega$, which coincide when the domain is the entire Euclidean space, but can otherwise be quite different. In particular, Cabre and Tan \cite{CT} have analyzed a very similar problem, taking as the fractional Laplacian the operator with the same eigenfunctions as the regular Laplacian, by extending to one further dimension. Another way is to restrict the integration to the domain:
\[ (-\Delta)^{\alpha/2} u(x) = C_{n,\alpha}  P.V. \int_{\Omega} \frac{u(x)-u(z)}{|x-z|^{n+\alpha}} dz , \]
known as the regional fractional Laplacian \cite{Gu}. In our paper, we mainly consider the cases when $\Omega$ is a unit ball or a half space, while our operator is defined by (\ref{Ad7});
or equivalently, by the Fourier transform: 
$$ \widehat{(-\lap)^{\alpha/2} u}(\xi) = |\xi|^{\alpha} \hat{u} (\xi) $$
where $\hat{u}$ is the Fourier transform of $u$. Obviously, this operator is well defined in
${\cal{S}}$, the Schwartz space of rapidly decreasing $C^{\infty}$ functions in $R^n$. One can
extend this definition to the distributions in the space
$$ {\cal{L}}_{\alpha/2} = \{ u \mid \int_{\mathbb{R}^n} \frac{|u(x)|}{1 + |x|^{n+\alpha}} d x < \infty \} $$
by
$$ < (-\lap)^{\alpha/2} u, \phi > = \int_{\mathbb{R}^n} u \, (-\lap)^{\alpha/2} \phi \, d x , \;\; \mbox{ for all } \; \phi \in C_0^{\infty}(\Omega) .$$

Given any $f \in L_{loc}^1(\Omega)$, we say that $u \in {\cal{L}}_{\alpha/2}$ solves the problem
$$ (-\lap)^{\alpha/2} u = f(x) , \;\; x \in \Omega $$
if and only if
\begin{equation}
\int_{\mathbb{R}^n} u \, (-\lap)^{\alpha/2} \phi \, d x = \int_{\mathbb{R}^n} f(x) \phi(x) d x , \;\; \mbox{ for all } \; \phi \in C_0^{\infty}(\Omega) .
\label{DS}
\end{equation}

Throughout this paper, we will consider the distributional solutions in the sense of (\ref{DS}).

Let $\textbf{B}_{1}=B_{1}(0)=\{x\in R^{n}:|x|<1\}$ be the unit ball in $\mathbb{R}^{n}$ and $0<\alpha<2$, $n\geq 3$. We first study the Dirichlet problem
\begin{equation}
\left\{\begin{array}{l}(-\Delta)^{\alpha/2}u=f(u),\;\;\;\;\;\mbox{in}\;\;\textbf{B}_{1},\\
u=0, \;\;\;\;\;\;\;\;
\;\;\;\;\;\;\;\;\;\;\;\;\;\;\;\;\mbox{in}\;\;\;\;\textbf{B}^{c}_{1}.
\end{array}
\right.
\label{1}
\end{equation}

Under some mild conditions on $f(\cdot)$, we will show that the positive solutions are radially symmetric and monotone
decreasing about the origin. Instead of using the extension method of Caffarelli and Silvestre \cite{CaS}, we introduce a new and
direct approach by studying an equivalent integral equation. We then use the {\em method of moving planes in integral forms} to
prove the symmetry and monotonicity of solutions.

Thanks to Kulczycki \cite{Ku}, under very mild regularity assumptions, for instance, $f(u) \in L^q(\textbf{B}_1)$ for some $q>1$, we can express the solutions of (\ref{1}) as
\begin{equation}\label{011}
u(x)=\int_{\textbf{B}_{1}}G_{1}(x,y)f(u(y))dy,
\end{equation}
where $G_{1}(x,y)$ is Green's function satisfying
$$
\left\{\begin{array}{l}(-\Delta)^{\alpha/2}G_1(x,y)=\delta(x-y),\;\;\;\;\;\mbox{in}\;\;\textbf{B}_{1},\\
G_1(x,y)=0, \;\;\;\;\;\;\;\;
\;\;\;\;\;\;\;\;\;\;\;\;\;\;\;\;\;\;\;\;\;\;\mbox{in}\;\;\;\;\textbf{B}^{c}_{1}.
\end{array}
\right.
$$
Set
$$s = |x-y|^2 \;\; \mbox{ and } \;\; t = (1-|x|^2)(1-|y|^2).$$
Then we can write the Green's function in the form
\begin{equation}
G_{1}(x,y)=\frac{A_{n,\alpha}}{s^{(n-\alpha)/2}}\left[1-\frac{B_{n,\alpha}}{(s+t)^{(n-2)/2}} \int_0^{\frac{s}{t}}
\frac{(s-tb)^{(n-2)/2}}{b^{\alpha/2}(1+b)} d b \right],\;\;x,y\in \textbf{B}_{1}.
\label{green}
\end{equation}
where $A_{n,\alpha}$ and $B_{n,\alpha}$ are constants depending on $n$ and $\alpha$.
\medskip

Our first result is
\begin{mthm}\label{mthm1}
Assume that  \\
$(f_{1})$ $f:[0,\infty)\rightarrow [0, \infty)$ is increasing, $f(0)=0,$ and either one of the following is satisfied: \\
$(f_{2})$ $f'(\cdot)$ is monotonic and $f'(u)\in L^{\frac{n}{\alpha}}(\textbf{B}_{1}),$
or\\
$(\widetilde{f}_{2})$
$$
|f'(u)|\leq C_{1}|u|^{\beta_{1}}+C_{2}|u|^{\beta_{2}}+C_{3},
$$
where $C_{1}$, $C_{2}$, $C_{3}$, and $\beta_{1}$ can be any nonnegative constants, while $\beta_{2}$
is some non-positive constant. If $C_{1}>0$, we require $|u|^{\beta_{1}}\in L^{\frac{n}{\alpha}}(\textbf{B}_{1})$, and if $C_{2}>0$, we need $|u|^{\beta_{2}}\in L^{\frac{n}{\alpha}}(\textbf{B}_{1})$.

Then every positive solution of integral equation (\ref{011}) is radially symmetric about the origin and strictly decreasing in the radial direction.
\end{mthm}

\begin{mcor}
Under the conditions of Theorem \ref{mthm1}, if $u$ is a positive solution of (\ref{1}) with $f(u) \in L^q(\textbf{B}_1)$ for some $q>1$, then it is radially symmetric about the origin and strictly decreasing in the radial direction.
\label{mcorH}
\end{mcor}

\begin{mrem}
i) Our prototype is $f(u) = u^p$. If $1<p\leq \frac{n+\alpha}{n-\alpha}$ and $u \in L^{\frac{2n}{n-\alpha}}(\textbf{B}_1)$ ( this can derived from $u \in H^{\alpha/2}$),
then all the conditions in the Corollary are met.

ii) If  $f(\cdot)$ is Lipschitz continuous, then it satisfies condition $(\tilde{f}_2)$ with $C_1=C_2=0.$

iii) Servadei and Valdinoci \cite{SV} in a recent paper obtained the classic existence result of Brezis and Nirenberg, demonstrating our class is non-empty. Fall and Weth \cite{FW} adapted the Pohozaev estimates to show nonexistence of solutions to the critical power equation itself ( see also \cite{RS}).

\end{mrem}

We also establish some regularity of the solutions.

\begin{mthm}\label{th2}
Let $u(x)$ be a positive solution of (\ref{1}) or of (\ref{011}). Assume that
\begin{equation} 
|\frac{f(u)}{u}| \leq C_1 + C_2 |u|^{\beta} , 
\label{R1}
\end{equation} 
for some $\beta > \frac{\alpha}{n-\alpha}$, and $u(x)\in L^{n\beta/\alpha}(\textbf{B}_{1})$.

Then $u$ is uniformly bounded in $\textbf{B}_{1}$.
\end{mthm}

\begin{mrem} Again, if $f(u) = u^p$ for $1<p\leq \frac{n+\alpha}{n-\alpha}$ with $u \in L^{\frac{2n}{n-\alpha}}(\textbf{B}_1)$, then one can verify that all the conditions in Theorem \ref{th2} are 
satisfied. 
\end{mrem} 

By applying some estimates found in \cite{CSS} and \cite{Si}, we apply a bootstrapping argument to show that $u$ is, in some sense, as smooth as $f$. In particular, if $f \in C^m(\mathbb{R})$, then $u \in C^m(\textbf{B}_{1})$.
\begin{mthm}
Assume that $u(x)$ is a positive bounded solution of (\ref{1}).  If $f:\mathbb{R}_+ \rightarrow \mathbb{R}$ is in $C^m$, then $u \in C^{m}(\textbf{B}_1)$.
\label{mthm2}
\end{mthm}
\bigskip

We then consider the case when $\Omega$ is the half space

\begin{equation}
\left\{\begin{array}{l}(-\Delta)^{\alpha/2}u(x)=u^{p}(x),\;\;\;\;\;\;\;\;\;\;\;\;\;\mbox{in}\;\;\mathbb{R}^n_+,\\
u(x)=0, \;\;\;\;\;\;\;\;
\;\;\;\;\;\;\;\;\;\;\;\;\;\;\;\;\;\;\;\;\mbox{in}\;\;\;\; \mathbb{R}_{-}^{n} ,
\end{array}
\right.
\label{hpde}
\end{equation}
where $\mathbb{R}^n_+ = \{ x = (x_1, x_2, \cdots, x_n) \in \mathbb{R}^n \mid x_n > 0\} ,$ and $\mathbb{R}_{-}^{n}$ is the complement of $\mathbb{R}^n_+$.

In order to use the method of integral equations, we first establish the equivalence between problem (\ref{hpde}) and the integral equation
\begin{equation}
u(x)=\int_{\mathbb{R}^n_+}G_{\infty}(x,y)u^{p}(y)dy,
\label{int1}
\end{equation}
where
\begin{equation}\label{04}
G_{\infty}(x,y)=\frac{A_{n,\alpha}}{s^{(n-\alpha)/2}}\left[1- \frac{B_{n,\alpha}}{(s+t)^{(n-2)/2}} \int_0^{\frac{s}{t}}
\frac{(s-tb)^{(n-2)/2}}{b^{\alpha/2}(1+b)} d b \right],\;\;x,y\in \mathbb{R}^n_+,
\end{equation}
is the Green function in $\mathbb{R}^n_+$ with the same Dirichlet condition. Here
$$s = |x-y|^2 \;\; \mbox{ while } \;\; t = 4 x_n y_n .$$
As compared to the Green's function $G_1(x,y)$ in the unit ball, the expression for $G_{\infty}(x,y)$ in terms of $s$ and $t$ is the same. However, $t$ is differently defined here.

We prove
\begin{mthm}
Assume that $u$ is a locally bounded positive solution of problem (\ref{hpde})
and the growth of $u$ is not as fast as a constant multiple of $(x_n)^{\alpha/2}$, more precisely, there exists a sequence $\{x^k\} \in \mathbb{R}_+^n$, such that
\begin{equation}
\frac{u(x^k)}{(x^k_n)^{\alpha/2}} \ra 0
\label{Ad40}
\end{equation}

Then $u$ is also a solution of integral equation (\ref{int1})
and vice versa.
\label{mthmA}
\end{mthm}

The proof of this theorem is based on the following uniqueness result.

\begin{mlem}

Assume that $w$ is a nonnegative solution of
\begin{equation}
\left\{
  \begin{array}{ll}
   (-\Delta)^{\alpha/2}w=0,\;\;& x\in \mathbb{R}_{+}^{n},\\
w= 0, &x\in \mathbb{R}_{-}^{n}.
  \end{array}
\right.
\label{Ad60}
\end{equation}

Then we have either
\begin{equation}
w(x)\equiv 0,\;x\in \mathbb{R}^{n},
\end{equation}
or there exists a constant $a_o > 0$, such that
\begin{equation}
w(x) \geq a_o (x_n)^{\alpha/2}, \;\; \forall \; x\in \mathbb{R}_+^{n}.
\label{Ad20}
\end{equation}
\end{mlem}

\begin{mrem}
The above lemma presents the best possible uniqueness result for such a problem. To see this,
for any constant $a_o$, let
$$
g(x) =\left\{
            \begin{array}{ll}
             a_o (x_{n})^{\alpha/2} ,& x \in \mathbb{R}^n_+,\\
0,& x \in \mathbb{R}^n_-.
            \end{array}
          \right.
$$

Then it is well known that $g(x)$ is a non-zero solution of problem (\ref{Ad60}) ( see \cite{CRS}).
\end{mrem}

Next, we establish Liouville theorems for the integral equation.
\begin{mthm}\label{mthm3}
  Assume $p>\frac{n}{n-\alpha}$. If $u\in L^{\frac{n(p-1)}{\alpha}}(\mathbb{R}^n_+)$ is a non-negative solution of integral equation (\ref{int1}), then $u(x)\equiv 0$.
\end{mthm}

Then by Theorem \ref{mthmA}, one derive immediately the following

\begin{mcor}
  Assume $p>\frac{n}{n-\alpha}$. If $u\in L^{\frac{n(p-1)}{\alpha}}(\mathbb{R}^n_+)$ is a non-negative solution of problem (\ref{hpde}), then $u(x)\equiv 0$.
\end{mcor}

To prove the non-existence of positive solutions for the integral equation, we again employ the {\em method of moving planes in integral forms}. We move the plane along $x_n$ direction to show that the solutions must be monotone increasing in $x_n$ and thus derive a contradiction.

To remarkably weaken the global integrability condition $u\in L^{\frac{n(p-1)}{\alpha}}(\mathbb{R}^n_+)$ in Theorem \ref{mthm3}, we exploit a
Kelvin type transform. To ensure that the half space $\mathbb{R}^n_+$ is invariant under the inversion, we need to place the centers on the boundary, $\partial \mathbb{R}^n_+$. For a point $z^{0}\in \partial \mathbb{R}^n_+$, we consider $$\bar{u}_{z^0}(x)=\frac{1}{|x-z^{0}|^{n-\alpha}}u \left(\frac{x-z^{0}}{|x-z^{0}|^{2}}+z^{0}\right),$$
the Kelvin type transform of $u(x)$ centered at $z^0$.  Some new ideas are involved.

In the critical case $p = \frac{n+\alpha}{n-\alpha}$, we consider two possibilities.

 (i) \emph{There is a point $z^{0}\in \partial \mathbb{R}^n_+$, such that $\bar{u}_{z^0}(x)$ is bounded near $z^{0}$}. In this situation, $u$ is globally integrable, and we move the planes in the direction of $x_n$-axis to show that the solution $u$ is monotone increasing in $x_n$, as we did in the proof of Theorem \ref{mthm3}.

 (ii) \emph{For all $z^{0}\in \partial \mathbb{R}^n_+$,  $\bar{u}_{z^0}(x)$ are unbounded near $z^{0}$.} In this situation, we move the planes in $x_1, \cdots, x_{n-1}$ directions to show that, for every $z^0$,  $\bar{u}_{z^0}$ is axially symmetric about the line that is parallel to $x_n$-axis and
 passing through $z^0$. This implies further that $u$ depends on $x_n$ only.

 In the subcritical case, we only need to work on $\bar{u}_{z^0}(x)$; and similar to the above possibility (ii), we show that for every $z^0$,  $\bar{u}_{z^0}$ is axially symmetric about the line that is parallel to $x_n$-axis and
 passing through $z^0$, which implies again that $u$ depends on $x_n$ only.

 In both cases, we will be able to derive a contradiction and prove the following

\begin{mthm}\label{mthm4}
Assume that $1 < p \leq\frac{n+\alpha}{n-\alpha}$. If $u$ is a non-negative locally bounded solution of (\ref{int1}), then $u(x)\equiv 0.$ In particular, when $p = \frac{n+\alpha}{n-\alpha}$, we only require $u \in L_{loc}^{\frac{2n}{n-\alpha}}(\mathbb{R}_+^n)$. 
\end{mthm}

\begin{mcor} Assume $1 < p \leq\frac{n+\alpha}{n-\alpha}$. If $u$ is a non-negative locally bounded solution of problem (\ref{hpde}) with the growth condition (\ref{Ad40}), then $u \equiv 0$.
\end{mcor}

\begin{mrem}
In \cite{FW}, to establish the non-existence of positive solutions for (\ref{hpde}) via the extension method, they required that $u \in \mathcal{D}^{\alpha/2, 2} \cap C(\mathbb{R}^n)$.
One can see that our growth condition here is much weaker.
\end{mrem}

It is well-known that these kinds of Liouville theorems play an important role in establishing a priori estimates for
the solutions of a family of corresponding boundary value problems in either bounded domains or Riemannian manifolds with
boundaries.

In Section 2, we study symmetry and monotonicity of positive solutions in the unit ball and establish Theorem \ref{mthm1}.
In Section 3, we obtain the regularity of solutions and prove Theorem \ref{th2} and \ref{mthm2}. In Section 4, we show the equivalence between
problem (\ref{hpde}) and integral equation (\ref{int1}), and in Section 5, we prove non-existence of positive solutions in the
half space $\mathbb{R}^n_+$ and thus establish Theorem \ref{mthm3} and \ref{mthm4}.

For more articles concerning the origin and applications of {\em the method of moving planes in integral forms}, please see \cite{CL3, CL4, CLO, CLO1, CLO2, CZ, FC, MCL} and the references therein.

\section{Symmetry of Solutions in the Ball}

In this section, we will use the {\em method of moving planes in integral forms } to obtain the radial symmetry and monotonicity of positive solution for integral equation (\ref{011}) and thus prove Theorem \ref{mthm1}, which leads to Corollary \ref{mcorH} immediately.  

\subsection{Properties of the Green's Function}

Let $\lambda\in (-1,0)$ be a real number and
$T_{\lambda}=\{x\in \mathbb{R}^n|x_{1}=\lambda\}$. We denote
$\Sigma_{\lambda}$ the region in the ball between the plane $x_{1}=-1$ and the
plane $x_{1}=\lambda$:
\begin{equation}\nonumber
\Sigma_{\lambda}=\{x=(x_{1},x_{2},\cdots,x_{n})\in
\textbf{B}_{1}|-1<x_{1}<\lambda\}.
\end{equation}
Let
\begin{equation}\nonumber
x^{\lambda}=(2\lambda-x_{1},x_{2},\cdots, x_{n})
\end{equation}
be the reflection of the point
$x=(x_{1},x_{2},\cdots,x_{n})$ about the plane
$T_{\lambda}$,
denote
$\Sigma_{\lambda}^{C}=\textbf{B}_{1}\backslash\Sigma_{\lambda}$, the complement of $\Sigma_{\lambda} $ in $\textbf{B}_{1}$,
and set $$
u_{\lambda}(x)=u(x^{\lambda})\;\;\mbox{and}\;\;w_{\lambda}(x)=u_{\lambda}(x)-u(x).
$$

To carry on the method of moving planes, we need the following
properties of the Green's function.
\begin{lem}\label{lemm1}
(i) For any $x,y\in \Sigma_{\lambda}$, $x\neq y$, we have
\begin{equation}\label{6}
G_{1}(x^{\lambda},y^{\lambda})>\max\{G_{1}(x^{\lambda},y),
G_{1}(x,y^{\lambda})\}
\end{equation}
and
\begin{equation}\label{7}
G_{1}(x^{\lambda},y^{\lambda})-G_{1}(x,y)>|G_{1}(x^{\lambda},y)-G_{1}(x,y^{\lambda})|.
\end{equation}
(ii) For any $x\in \Sigma_{\lambda}$, $y\in \Sigma_{\lambda}^{C}$, it holds
\begin{equation}\label{91}
G_{1}(x^{\lambda},y)> G_{1}(x,y).
\end{equation}
\end{lem}
\textbf{Proof.}
$(i)$ Let $s= s(x,y)=|x-y|^{2}$. Since $x,y\in\Sigma_{\lambda}$, it is easy to verify that
\begin{equation}\label{8}
s(x^{\lambda},y^{\lambda})<s(x,y^{\lambda}),
\end{equation}
\begin{equation}\label{0.12}
s(x,y)<s(x,y^{\lambda}),
\end{equation}
\begin{equation}\label{9}
s(x^{\lambda},y)=s(x,y^{\lambda}),
\end{equation}
\begin{equation}\label{0.1}
s(x,y)=s(x^{\lambda},y^{\lambda}).
\end{equation}

Set $t = t(x,y)=(1-|x|^{2})(1-|y|^{2})$. Obviously, we have
\begin{eqnarray}\label{010}
t(x^{\lambda},y^{\lambda})&>&\max\{t(x,y^{\lambda}),t(x^{\lambda},y)\}\\
&\geq&\min\{t(x,y^{\lambda}),t(x^{\lambda},y)\}\\
\label{0.13}
&>&t(x,y).
\end{eqnarray}

Let  $B_{n, \alpha} = B$,
$$
I_{1}(s,t)=\frac{1}{(t+s)^{\frac{n-2}{2}}}
\int_{0}^{\frac{s}{t}}\frac{(s-tb)^{\frac{n-2}{2}}}{b^{\alpha/2}(1+b)}db,
$$
and
$$
H(s,t)=\frac{1}{s^{\frac{n-\alpha}{2}}}\left[1-B I_{1}(s,t)\right].
$$
Then we can rewrite
$$
G_{1}(x,y) = A_{n,\alpha}H(s,t).
$$

Since
\begin{equation}
\label{1aa}
\frac{\partial I_{1}}{\partial s} = \frac{(n-2)t}{2(s+t)^{\frac{n}{2}}}
\int_{0}^{\frac{s}{t}}\frac{(s-tb)^{\frac{n-4}{2}}}{b^{\alpha/2}}db>0,
\end{equation}
we have
\begin{equation}\label{13}
\frac{\partial H}{\partial s}=-\frac{n-\alpha}{2s^{\frac{n-\alpha}{2}+1}}\left[1-BI_{1}(s,t)\right]+\frac{-B}{s^{\frac{n-\alpha}{2}}}\frac{\partial I_{1}}{\partial s}<0.
\end{equation}
By a straight forward calculation,
\begin{equation}
\label{1ab}
\frac{\partial I_{1}}{\partial t}
=-\frac{(n-2)s}{2(t+s)^{\frac{n}{2}}}\int_{0}^{\frac{s}{t}}\frac{(s-tb)^{\frac{n-4}{2}}}{b^{\alpha/2}}db<0,
\end{equation}
and consequently,
\begin{equation}
\label{14z}
\frac{\partial H}{\partial t}=\frac{-B}{s^{\frac{n-\alpha}{2}}}\frac{\partial I_{1}}{\partial t}>0.
\end{equation}
Let $b=s\lambda$. we obtain
\begin{equation}
\frac{\partial I_{1}}{\partial t}=-\frac{(n-2)s^{\frac{n-\alpha}{2}}}{2(t+s)^{\frac{n}{2}}}\int_{0}^{\frac{1}{t}}
\frac{(1-t\lambda)^{\frac{n-4}{2}}}{\lambda^{\alpha/2}}d\lambda.
\end{equation}
Then we have
\begin{equation}\label{10}
\frac{\partial H}{\partial t}
=\frac{(n-2)B}{2(t+s)^{\frac{n}{2}}}\int_{0}^{\frac{1}{t}}
\frac{(1-t\lambda)^{\frac{n-4}{2}}}{\lambda^{\alpha/2}}d\lambda.
\end{equation}
It follows that
\begin{equation}\label{81}
\frac{\partial^{2}H}{\partial t\partial s}<0.
\end{equation}
Since
$$
G_{1}(x,y)=A_{n,\alpha}H(s(x,y), t (x,y)),
$$
we derive (\ref{6}) from (\ref{8})-(\ref{010}), (\ref{13}), and  (\ref{14z}).

Next, we prove (\ref{7}). It can be obtained from (\ref{0.12})-(\ref{9}),
(\ref{010})-(\ref{0.13}), and (\ref{81}) as follows:
\begin{eqnarray}\nonumber
&&G_{1}(x^{\lambda},y^{\lambda})-G_{1}(x,y)\\
\nonumber
&=&A_{n,\alpha}\int_{t(x,y)}^{t(x^{\lambda},y^{\lambda})}\frac{\partial
H(s(x,y),t)}{\partial t}dt\\
\nonumber
&>&A_{n,\alpha}\int_{t(x,y)}^{t(x^{\lambda},y^{\lambda})}\frac{\partial
H(s(x^{\lambda},y),t)}{\partial t}dt\\
\nonumber &\geq&
A_{n,\alpha}\int_{t(x,y^{\lambda})}^{t(x^{\lambda},y)}\frac{\partial
H(s(x^{\lambda},y),t)}{\partial t}dt\\
\nonumber
&=&A_{n,\alpha}[ H(s(x^{\lambda},y),t(x^{\lambda},y))-H(s(x,y^{\lambda}),t(x,y^{\lambda}))]\\
&=&G_{1}(x^{\lambda},y)-G_{1}(x,y^{\lambda}).
\end{eqnarray}

Similarly,
\begin{eqnarray} \nonumber
&&G_1(x^{\lambda}, y^{\lambda}) - G_1(x,y) \\
\nonumber
&>&A_{n,\alpha}\int_{t(x,y)}^{t(x^{\lambda},y^{\lambda})}\frac{\partial
H(s(x,y^{\lambda}),t)}{\partial t}dt\\
\nonumber &\geq&
A_{n,\alpha}\int_{t(x^{\lambda},y)}^{t(x,y^{\lambda})}\frac{\partial
H(s(x, y^{\lambda}),t)}{\partial t}dt\\
\nonumber
&=&A_{n,\alpha}[ H(s(x,y^{\lambda}),t(x,y^{\lambda}))-H(s(x^{\lambda},y),t(x^{\lambda},y))]\\
\nonumber
&=& G_1(x, y^{\lambda}) - G_1(x^{\lambda},y) .
\end{eqnarray}

$(ii)$ For $x\in \Sigma_{\lambda}$ and $y\in \Sigma_{\lambda}^{C}$, we have
\begin{equation}\label{52}
s(x^{\lambda},y)<s(x,y),\;\;t(x^{\lambda},y)>t(x,y).
\end{equation}
Then by (\ref{13})-(\ref{14z}), and (\ref{52}), we deduce
$$
G_{1}(x^{\lambda},y)> G_{1}(x,y).
$$
This completes the proof of Lemma \ref{lemm1}.

\begin{lem}\label{lemm2}
For any $x\in \Sigma_{\lambda}$, it holds
\begin{equation}
u(x)-u_{\lambda}(x)\leq \int_{\Sigma_{\lambda}}[G_{1}(x^{\lambda},y^{\lambda})-G_{1}(x,y^{\lambda})][f(u(y))-f(u_{\lambda}(y))]dy.
\end{equation}
\end{lem}
\textbf{Proof.} Since
$$ G(x,y) f(u) \mid_{\tilde{\Sigma}_{\lambda}} = G(x, y^{\lambda}) f(u_{\lambda}) \mid_{\Sigma_{\lambda}},$$
we have
\begin{eqnarray*}
u(x)&=&\int_{\Sigma_{\lambda}}G_{1}(x,y)f(u(y))dy  \\
    &+& \int_{\Sigma_{\lambda}}G_{1}(x,y^{\lambda})f(u_{\lambda}(y))dy
+\int_{\Sigma_{\lambda}^{C}\backslash\widetilde{\Sigma}_{\lambda}}G_{1}(x,y)f(u(y))dy,
\end{eqnarray*}
and
\begin{eqnarray*}
u(x^{\lambda})&=& \int_{\Sigma_{\lambda}}G_{1}(x^{\lambda},y)f(u(y))dy \\
&+& \int_{\Sigma_{\lambda}}G_{1}(x^{\lambda},y^{\lambda})f(u_{\lambda}(y))dy+
\int_{\Sigma_{\lambda}^{C}\backslash\widetilde{\Sigma}_{\lambda}}G_{1}(x^{\lambda},y)f(u(y))dy,
\end{eqnarray*}
where $\widetilde{\Sigma}_{\lambda}=\{x^{\lambda} \mid x\in \Sigma_{\lambda}\}$ is the reflection of $\Sigma_{\lambda}$ about the plane $T_{\lambda}$.
By Lemma \ref{lemm1}, we arrive at
\begin{eqnarray}\nonumber
u(x)-u(x^{\lambda})&=&\int_{\Sigma_{\lambda}}\left[G_{1}(x,y)-G_{1}(x^{\lambda},y)\right]f(u(y))dy\\
\nonumber
&&+\int_{\Sigma_{\lambda}}\left[G_{1}(x,y^{\lambda})-G_{1}(x^{\lambda},y^{\lambda})\right]f(u_{\lambda}(y))dy\\
\nonumber
&&+\int_{\Sigma_{\lambda}^{C}\backslash\widetilde{\Sigma}_{\lambda}}\left[G_{1}(x,y)-G_{1}(x^{\lambda},y)\right]f(u(y))dy\\
\nonumber
&\leq&\int_{\Sigma_{\lambda}}\left[G_{1}(x,y)-G_{1}(x^{\lambda},y)\right]f(u(y))dy\\
\nonumber
&&-\int_{\Sigma_{\lambda}}\left[G_{1}(x^{\lambda},y^{\lambda})-G_{1}(x,y^{\lambda})\right]f(u_{\lambda}(y))dy\\
\nonumber
&\leq&\int_{\Sigma_{\lambda}}\left[G_{1}(x^{\lambda},y^{\lambda})-G_{1}(x,y^{\lambda})\right]f(u(y))dy\\
\nonumber
&&-\int_{\Sigma_{\lambda}}\left[G_{1}(x^{\lambda},y^{\lambda})-G_{1}(x,y^{\lambda})\right]f(u_{\lambda}(y))dy\\
\nonumber &=&\int_{\Sigma_{\lambda}}\left[
G_{1}(x^{\lambda},y^{\lambda})-G_{1}(x,y^{\lambda})\right]\left[f(u(y))-f(u_{\lambda}(y))\right]dy.
\end{eqnarray}

This completes the proof of the lemma.

\begin{lem} \label{le1} (An equivalent form of the Hardy-Littlewood-Sobolev inequality)  Assume $0<\alpha<n$ and $\Omega\subset \mathbb{R}^n$. Let $g\in L^{\frac{np}{n+\alpha p}}(\Omega)$ for $\frac{n}{n-\alpha}<p<\infty.$ Define
$$
Tg(x):=\int_{\Omega}\frac{1}{|x-y|^{n-\alpha}}g(y)dy.
$$
Then
\begin{equation}\label{f1} \|Tg\|_{L^{p}(\Omega)}\leq C(n,p, \alpha)\|g\|_{L^{\frac{np}{n+\alpha p}}(\Omega)}.
\end{equation} \end{lem}
The proof of this lemma is standard and can be found in \cite{CL1} or \cite{CL2}.

\subsection{The Proof of Theorem \ref{mthm1}}

The proof of Theorem \ref{mthm1} consists of two steps. In \emph{Step 1.}, we show that for $\lambda$
sufficiently close to $-1$, we have
\begin{equation}\label{2a}
w_{\lambda}(x)=u_{\lambda}(x)-u(x)\geq 0, \;\;a.e.\;\;x\in \Sigma_{\lambda}.
\end{equation}
This provides a starting point to move the plane $T_{\lambda}$ along the $x_{1}$ direction.
In \emph{Step 2.}, we move the plane continuously to the right as long as inequality (\ref{2a}) holds. We
show that the plane can be moved all the way to $\lambda=0$ and thus derive
\begin{equation}\label{2b}
u(-x_{1},x_{2},\cdots,x_{n})\leq u(x_{1},x_{2},\cdots,x_{n}),\;\;\forall x\in \textbf{B}_{1},\;\;x_{1}\geq 0.
\end{equation}

Similarly, we can start the plane $T_{\lambda}$ from close to $\lambda=1$ and move it to the left to the limiting
position $T_{0}$ to deduce
\begin{equation}\label{2c}
u(-x_{1},x_{2},\cdots,x_{n})\geq u(x_{1},x_{2},\cdots,x_{n}),\;\;\forall x\in \textbf{B}_{1},\;\;x_{1}\geq 0.
\end{equation}
Now (\ref{2b}) and (\ref{2c}) imply that $u(x)$ is symmetric about the plane $T_{0}$. Since the direction
of $x_{1}$ can be chosen arbitrary, we deduce that $u(x)$ is radially symmetric about the origin and strictly decreasing in the
radial direction.

\emph{Step 1.} Define
$$
\Sigma_{\lambda}^{-}=\{x\in\Sigma_{\lambda}| \;\;u(x)>u_{\lambda}(x)\}.
$$
We show that $\Sigma_{\lambda}^{-}$ is
almost empty by estimating a certain integral norm on it.

By Lemma \ref{lemm1}, \ref{lemm2}, Mean Value Theorem, and $(f_{1})$, we have, for any $x\in\Sigma_{\lambda}^{-}$,
\begin{eqnarray}
\nonumber
0&<&u(x)-u_{\lambda}(x)\\
\nonumber
&\leq& \int_{\Sigma_{\lambda}}[G_{1}(x^{\lambda},y^{\lambda})-G_{1}(x,y^{\lambda})][f(u(y))-f(u_{\lambda}(y))]dy\\
\nonumber
&\leq& \int_{\Sigma_{\lambda}^{-}}[G_{1}(x^{\lambda},y^{\lambda})-G_{1}(x,y^{\lambda})][f(u(y))-f(u_{\lambda}(y))]dy\\
\nonumber
&\leq& \int_{\Sigma_{\lambda}^{-}}G_{1}(x^{\lambda},y^{\lambda})[f(u(y))-f(u_{\lambda}(y))]dy\\
\label{2d}
&\leq& \int_{\Sigma_{\lambda}^{-}}G_{1}(x^{\lambda},y^{\lambda})|f'(\psi(y))||w_{\lambda}(y)|dy
\end{eqnarray}
where $\psi(y)$ is valued between $u(y)$ and $u_{\lambda}(y)$. By the formula for $G_{1}(x,y)$, it is easy to see
\begin{equation}\label{2e}
|G_{1}(x,y)|\leq\frac{C}{|x-y|^{n-\alpha}}.
\end{equation}
It follows from (\ref{2d}) and (\ref{2e}) that, for any $x\in \Sigma_{\lambda}^{-}$,
\begin{equation}
0<u(x)-u_{\lambda}(x)\leq C\int_{\Sigma_{\lambda}^{-}}\frac{1}{|x-y|^{n-\alpha}}|f'(\psi(y))||w_{\lambda}(y)|dy.
\end{equation}
Applying the HLS inequality (Lemma \ref{le1}) and H\"{o}lder inequality,
we have, for any $q>\frac{n}{n-\alpha}$,
\begin{eqnarray}
\nonumber
\|w_{\lambda}\|_{L^{q}(\Sigma_{\lambda}^{-})}&\leq& C\|f'(\psi(x))w_{\lambda}\|_{L^{\frac{nq}{n+\alpha q}(\Sigma_{\lambda}^{-})}}\\
\label{ab}
&\leq & C\|f'(\psi(x))\|_{L^{\frac{n}{\alpha}}(\Sigma_{\lambda}^{-})}\|w_{\lambda}(x)\|_{L^{q}(\Sigma_{\lambda}^{-})}.
\end{eqnarray}
By assumption $(f_{2})$ or $(\widetilde{f}_{2})$, for $\lambda$ sufficiently close to $-1$, we have
$$
C\|f'(\psi(x))\|_{L^{\frac{n}{\alpha}}(\Sigma_{\lambda}^{-})}\leq\frac{1}{2}.
$$
Then by (\ref{ab}), we have
$$
\|w_{\lambda}\|_{L^{q}(\Sigma_{\lambda}^{-})}\leq\frac{1}{2}\|w_{\lambda}\|_{L^{q}(\Sigma_{\lambda}^{-})}.
$$
This implies that $\|w_{\lambda}\|_{L^{q}(\Sigma_{\lambda}^{-})}=0$, therefore $\Sigma_{\lambda}^{-}$ must
be measure zero. And then (\ref{2a}) holds.

\emph{Step 2.} We now move the plane $T_{\lambda}$ continuously towards the right as long
as inequality (\ref{2a}) holds to its limiting position.

Define
$$
\lambda_{0}=\sup\{\lambda\in (-1,0)|\;\;w_{\kappa}(x) \geq 0,\;\;x\in\Sigma_{\kappa},\;\;\kappa\leq\lambda\}.
$$
We prove that $\lambda_{0}$ must be $0$.

Otherwise, suppose $\lambda_{0}<0$. First, we shall show that
\begin{equation}\label{ab25}
w_{\lambda_{0}}(x)>0
\end{equation}
in the interior of $\Sigma_{\lambda_{0}}$.

Indeed, by the first two expressions in the proof of Lemma \ref{lemm2} and Lemma \ref{lemm1}, we have
\begin{eqnarray}
\nonumber
u_{\lambda}(x)-u(x)&\geq&\int_{\Sigma_{\lambda}}\left[G_{1}(x^{\lambda},y^{\lambda})
-G_{1}(x,y^{\lambda})\right]\left[f(u_{\lambda}(y))-f(u(y))\right]dy\\
\nonumber
&\;\;\;\;+&\int_{\Sigma_{\lambda}^{C}\backslash\widetilde{\Sigma}_{\lambda}}
\left[G_{1}(x^{\lambda},y)-G_{1}(x,y)\right]f(u(y))dy\\
\label{m1}
&\geq& \int_{\Sigma_{\lambda}^{C}\backslash\widetilde{\Sigma}_{\lambda}}
\left[G_{1}(x^{\lambda},y)-G_{1}(x,y)\right]f(u(y))dy.
\end{eqnarray}
If (\ref{ab25}) is violated, then there exists some point
$x_{0}\in\Sigma_{\lambda_{0}}$ such that
$u(x_{0})=u_{_{\lambda_{0}}}(x_{0})$. Consequently, by (\ref{91}) and (\ref{m1}), we obtain
\begin{equation}\label{ab26}
f(u(y))\equiv 0,\;\;\forall y\in\Sigma_{\lambda_{0}}^{C}\backslash\widetilde{\Sigma}_{\lambda_{0}} .
\end{equation}
Due to $(f_{1})$, we must have
$$
u\equiv 0,\;\;\forall \; y\in\Sigma_{\lambda_{0}}^{C}\backslash\widetilde{\Sigma}_{\lambda_{0}}.
$$
This is a contradiction with our assumption that $u>0$. Therefore (\ref{ab25}) must be true.

By virtue of the Lusin Theorem, for any $\delta>0$, there exists a
closed subset $F_{\delta}$ of $\Sigma_{\lambda_{0}}$,
with $\mu(\Sigma_{\lambda_{0}}\backslash F_{\delta})<\delta$, such
that $w_{\lambda_{0}}|_{F_{\delta}}$ is continuous with respect to $x$, and
hence $w_{\lambda}|_{F_{\delta}}$ is continuous with respect to $\lambda$
for $\lambda$ close to $\lambda_{0}$. By (\ref{ab25}), there
exists $\epsilon>0$ such that for all $\lambda\in [\lambda_{0},\lambda_{0}+\epsilon)$, it
holds
$$
w_{\lambda}(x)\geq 0,\;\;\forall x\in F_{\delta}.
$$
It follows that, for such $\lambda$,
$$
\mu(\Sigma_{\lambda}^{-})\leq \mu(\Sigma_{\lambda_{0}}\backslash F_{\delta})+\mu(\Sigma_{\lambda}\backslash\Sigma_{\lambda_{0}})\leq\delta +2\epsilon.
$$
Similar to \emph{Step 1}, we can choose $\delta$ and $\epsilon$ sufficiently small such that
$$
C\|f'(\psi(x))\|_{L^{\frac{n}{\alpha}}(\Sigma_{\lambda}^{-})}\leq\frac{1}{2}.
$$
Consequently from (\ref{ab}), we have $\|w_{\lambda}(x)\|_{L^{q}(\Sigma_{\lambda}^{-})}=0$, and
hence $\Sigma_{\lambda}^{-}$ must be measure zero, and hence
$$
w_{\lambda}(x)\geq 0,\;\;a.e.\;\;x\in \Sigma_{\lambda},\;\;\lambda\in [\lambda_{0},\lambda_{0}+\epsilon).
$$
This contradicts the definition of $\lambda_{0}$. Therefore we must have $\lambda_{0}=0$.
This completes the proof of Theorem \ref{mthm1}.

\begin{rem}

In the above, we only presented the proof under condition $(f_2)$, while a similar argument works under condition
$(\tilde{f}_2)$.

\end{rem}

\subsection{The Proof of Corollary \ref{mcorH}}

By Theorem \ref{mthm1}, we only need to show that under the regularity assumption
\begin{equation}
f(u) \in L^q (\textbf{B}_1) \; \mbox{ for some } q > 1 , 
\label{Z01}
\end{equation}
$u$ satisfies integral equation (\ref{011}). 

\emph{Proof}. Let 
$$ v(x) = \int_{\textbf{B}_1} G_1(x,y) f(u(y)) d y .$$ 
Then obviously, 
$$ |v(x)| \leq C \int_{\textbf{B}_1} \frac{1}{|x-y|^{n-\alpha}} |f(u(y))| d y . $$
Consequently, by HLS inequality, we have 
$$ \|v\|_{L^r(\textbf{B}_1)} \leq C \|f(u)\|_{L^q(\textbf{B}_1)} ,$$
where $ r = \frac{nq}{n-\alpha q} > 1$. 

Moreover, one can easily verify that 
$$
\left\{\begin{array}{ll} 
(-\lap)^{\alpha/2} v = f(u) , & \mbox{ in } \textbf{B}_1 , \\
v = 0 , & \mbox{ on } \mathbb{R}^n \setminus \textbf{B}_1 . 
\end{array}
\right. 
$$

Let $w = u-v$, then 
$$
\left\{\begin{array}{ll} 
(-\lap)^{\alpha/2} w = 0 , & \mbox{ in } \textbf{B}_1 , \\
w = 0 , & \mbox{ on } \mathbb{R}^n \setminus \textbf{B}_1 . 
\end{array}
\right. 
$$

Recall that $w \in L^1$ is fractional harmonic in the weak sense
\begin{equation} 
\int_{\mathbb{R}^n} w \, (- \lap)^{\alpha/2} \phi (x) d x = 0 , \;\; \forall \; \phi \in C_0^{\infty}(\textbf{B}_1).
\label{R10}
\end{equation} 
We show that $w \equiv 0$. Otherwise, there exists a smooth, compactly supported function $\psi$ such that 
$$ \int_{\textbf{B}_1} w \, \psi \, d x \geq \delta > 0 .$$

Let 
$$ \phi(x) = \int_{\textbf{B}_1} G_1(x,y) \psi (y) d y . $$
Then
$$ \int_{\mathbb{R}^n} w \, (- \lap)^{\alpha/2} \phi (x) d x = \int_{\textbf{B}_1} w \, \psi \, d x \geq \delta > 0 $$
which is a contradiction with (\ref{R10}). 
This completes the proof of the corollary.

\section{Regularity of Solutions}

In this section, we establish regularity for positive solutions of
(\ref{1}) and of (\ref{011}), in which  
the following lemma from \cite{CL2} is a key ingredient.

Let $V$ be a Hausdorff topological vector space. Suppose there are two
extended norms defined on $V$,
$$
\|\cdot\|_{X},\;\;\|\cdot\|_{Y}: V\rightarrow  [0,\infty].
$$
Assume that the spaces
$$
X:=\{v\in V:\|v\|_{X}<\infty\}\;\;\mbox{and}\;\;Y:=\{v\in V:\|v\|_{Y}<\infty\}
$$
are complete under the corresponding norms, and the convergence in $X$ or in $Y$
implies the convergence in $V$.

\begin{lem}\label{lem3}(Regularity Lifting)
Let $T$ be a contracting map from $X$ into itself and from $Y$ into itself. Assume that
$f\in X$, and that there exists a function $g\in Z:=X\bigcap Y$ such that
$f=Tf+g$ in $X$. Then $f$ also belongs to $Z$.
\end{lem}

We have
\begin{thm}\label{addth2}
Let $u(x)$ be a positive solution of (\ref{1}) or of (\ref{011}). Assume that
\begin{equation} 
|\frac{f(u)}{u}| \leq C_1 + C_2 |u|^{\beta} , 
\label{m3}
\end{equation} 
for some $\beta > \frac{\alpha}{n-\alpha}$, and 
\begin{equation}
u(x)\in L^{n\beta/\alpha}(\textbf{B}_{1}).
\label{m6}
\end{equation} 

Then $u$ is uniformly bounded in $\textbf{B}_{1}$.
\end{thm}

\emph{Proof.}
We first work on the solutions of integral equation (\ref{011}). We will use the {\em Regularity Lifting} Lemma 
to show that
\begin{equation}\label{3a}
u\in L^{p}(\textbf{B}_{1}),\;\;\mbox{for}\;\;\mbox{any}\;\;p>\frac{n}{n-\alpha}.
\end{equation}

For any real number $a>0$, let $A=\{x\in \textbf{B}_{1}|\;\;u(x)>a\}$ and
$$
u_{a}(x)=\left\{\begin{array}{l} u(x),\;\;\;\;\;\;\mbox{if}\;\;x\in A,\\
0,\;\;\;\;\;\;\;\;\;\;\;\mbox{elsewhere}.\;\;
\end{array}
\right.
$$
Set $u_{b}=u(x)-u_{a}(x)$. Then obviously
$$
f(u)=f(u_{a})\chi_{A}+f(u_{b})\chi_{D},
$$
where $\chi_{A}$ is the characteristic function on the set A and $D=\textbf{B}_{1}\backslash A$.

Define the linear operator
$$
T_{a}w(x)=\int_{\textbf{B}_{1}}G_{1}(x,y)\frac{f(u_{a}(y))\chi_{A}(y)}{u_{a}(y)}w(y)dy,
$$
$$
J(x)=\int_{\textbf{B}_{1}}G_{1}(x,y)f(u_{b}(y))\chi_{D}(y)dy.
$$
Then obviously, $u$ satisfies the equation
$$
u(x)=T_{a}u(x)+J(x),\;\;\forall x\in \textbf{B}_{1}.
$$
We prove that, for $a$ sufficiently large, $T_{a}$ is a contracting map from $L^{p}(\textbf{B}_{1})$ to
$L^{p}(\textbf{B}_{1})$, for any $p>\frac{n}{n-\alpha}$. In fact, by (\ref{2e}), (\ref{m3}), HLS inequality,  and H\"{o}lder inequality, we have
\begin{eqnarray*} 
\|T_{a}w\|_{L^{p}(\textbf{B}_{1})} &\leq& C\|\frac{f(u_{a})\chi_{A}}{u_{a}}w\|_{L^{\frac{np}{n+\alpha p}}(\textbf{B}_{1})} \\
&\leq& C \|\frac{f(u)}{u} \|_{L^{\frac{n}{\alpha}}(A)}\|w\|_{L^{p}(\textbf{B}_{1})} \\
&\leq& C \|u\|_{L^{n\beta/\alpha}(A)} \|w\|_{L^{p}(\textbf{B}_{1})}.
\end{eqnarray*}

By (\ref{m6}), we can choose $a$ sufficiently large, so that the measure of $A$ is small
and hence
$$
\|T_{a}w\|_{L^{p}(\textbf{B}_{1})}\leq \frac{1}{2}\|w\|_{L^{p}(\textbf{B}_{1})}.
$$
Therefore $T_{a}$ is a contracting operator from  $L^{p}(\textbf{B}_{1})$ to $L^{p}(\textbf{B}_{1})$.

To estimate $J(x)$, we apply HLS inequality, (\ref{m3}), and the fact that $u \leq a$ on $D$ to derive
$$
\|J(x)\|_{L^{p}(\textbf{B}_{1})} \leq C\|f(u_{b})\chi_{D}\|_{L^{\frac{np} 
{n+\alpha p}}(\textbf{B}_{1})}\leq C_a.
$$
Hence $J(x)\in L^{p}(\textbf{B}_{1})$ for any $p>\frac{n}{n-\alpha}$. 

From our assumption, $u$ is in $L^{n\beta/\alpha}(\textbf{B}_1)$ with $n\beta/\alpha > \frac{n}{n-\alpha}$.
Then by Lemma \ref{lem3}, we arrive at (\ref{3a}).

By H\"{o}lder inequality, (\ref{m3}), and (\ref{3a}), we have
\begin{eqnarray*}
\nonumber
u(x)&=&\int_{\textbf{B}_{1}}G_{1}(x,y)f(u(y))dy\\
&\leq C&\left(\int_{\textbf{B}_{1}}\frac{1}{|x-y|^{a(n-\alpha)}}dy\right)^{\frac{1}{a}}
\left(\int_{\textbf{B}_{1}} \left( C_1 u + C_2 u^{(\beta+1)} \right)^b dy\right)^{\frac{1}{b}}\\
&\leq& C \left(\int_{\textbf{B}_{1}}\frac{1}{|x-y|^{a(n-\alpha)}}dy\right)^{\frac{1}{a}}
\end{eqnarray*}
with $\frac{1}{a}+\frac{1}{b}=1$. Choosing $1<a<\frac{n}{n-\alpha}$, the above integral is uniformly
bounded for all $x \in \textbf{B}_1$. Therefore $u$ is uniformly bounded. 

Next, we show that this is also true for equation (\ref{1}). In fact, by (\ref{m3}) and (\ref{m6}), we have 
$$ f(u) \in L^{\frac{n\beta}{\alpha (\beta +1)}} (\textbf{B}_1) . $$
Since we assume that $\beta > \frac{\alpha}{n-\alpha}$, it is easy to veryfy that $\frac{n\beta}{\alpha (\beta +1)} > 1.$ Now it follows from the proof of Corollary \ref{mcorH} that if $u$ is a solution of (\ref{1}), then it is also a solution of (\ref{011}), therefore it is uniformly bounded. 

This completes the proof of Theorem \ref{addth2}.

\begin{thm}
Assume that $u$ is a positive bounded solution of (\ref{1}).  If $f:\mathbb{R}_+ \rightarrow \mathbb{R}$ is in $C^m$, then $u \in C^{m}(\textbf{B}_1)$.
\end{thm}
\emph{Proof.}
Let
\[ w(x) = (-\Delta)^{-\alpha/2} (f(u)\chi_{\textbf{B}_1}) = C_{n,\alpha} \int_{\mathbb{R}^n} \frac{f(u(y)) \chi_{\textbf{B}_1}(y) }{|x-y|^{n-\alpha}} dy, \]
the Riesz potential of $f(u)$. Since $f(u)$ is both bounded and compactly supported, we know that its Riesz potential is well defined and is in $C^{\gamma}$ for any $\gamma < \alpha$ \cite{Si}.
Let $v = u - w$. Then $(-\Delta)^{\alpha/2} v = 0$ for $ x \in \textbf{B}_{1}$, with $v = -w$ in $\textbf{B}_{1}^c$. Fractional harmonic functions are known to be $C^\infty$ in the domain where they are fractionally harmonic, since their derivatives satisfy the same equation (see, for example, corollary 2.5 of \cite{CSS}).

Hence, we have $u \in C^\gamma$ in any domain compactly contained within $\textbf{B}_{1}$. To iterate in general, suppose $u \in C^{k,\gamma}$ in any domain compactly contained in $\textbf{B}_{1}$. If $f$ is $C^m$ for some $m\geq k$, then $f(u)$ is $C^{k,\gamma}$ as well. To iterate our scheme, we apply the estimates from \cite{Si} which state that if $f(u) \in C^{k,\gamma}$, then its Riesz potential is in $C^{l,\beta}$, where $l$ is the integer part of $k + \alpha + \gamma$, and $\beta$ the fractional part. We repeat until we get that $u \in C^{m,\gamma}$ for some $\gamma$, in every domain compactly contained within $\textbf{B}_{1}$, which tells us that $u \in C^m(\textbf{B}_1)$.

\section{Equivalence between the Two Equations on $\mathbb{R}^n_+$}

We first derive the expression of the Green's function in this half space. In the previous sections, we introduced from \cite{Ku} the Green's function of the operator $(-\Delta)^{\alpha/2}$ with Dirichlet
conditions on the unit ball $\textbf{B}_{1}$:
$$
G_{1}(x,y)= \frac{A_{n,\alpha}}{s^{\frac{n-\alpha}{2}}}\left[1-
B\frac{1}{(t+s)^{\frac{n-2}{2}}}
\int_{0}^{\frac{s}{t}}\frac{(s-tb)^{\frac{n-2}{2}}}{b^{\alpha/2}(1+b)}db\right]
$$
where $s=|x-y|^{2}$ and $t=(1-|x|^{2})(1-|y|^{2})$.

Set $P_{R}:=(0,\cdots,R)\in \mathbb{R}^n_+$, and $B_{R}(P_{R}):=\{x\in \mathbb{R}^n:|x-P_{R}|<R\}$, the ball
of radius $R$ centered at $P_{R}$. Let
$$
s_{R}=|\frac{x-P_{R}}{R}-\frac{y-P_{R}}{R}|^{2}=\frac{|x-y|^{2}}{R^{2}}
$$
and
$$
t_{R}=\left(1-|\frac{x-P_{R}}{R}|^{2}\right)\left(1-|\frac{y-P_{R}}{R}|^{2}\right)=
\left(\frac{2x_{n}}{R}-\frac{|x|^{2}}{R^{2}}\right)\left(\frac{2y_{n}}{R}-\frac{|y|^{2}}{R^{2}}|\right),
$$

Then we can write the Green's function on $B_{R}(P_{R})$ as
\begin{eqnarray}
\nonumber
&&G_{R}(x,y)\\
\nonumber
&=&\frac{1}{R^{n-\alpha}}G_{1}(\frac{x-P_{R}}{R},\frac{y-P_{R}}{R})\\
\label{36}
&=&\frac{A_{n,\alpha}}{|x-y|^{n-\alpha}}\left[1-
B\int_{0}^{\frac{s_{R}}{t_{R}}}\frac{\left(\frac{1-\frac{t_{R}}{s_{R}}b}
{1+\frac{t_{R}}{s_{R}}}\right)^{\frac{n-2}{2}}}{b^{\alpha/2}(1+b)}db\right].
\end{eqnarray}

Let $R\rightarrow\infty$ in (\ref{36}), we arrive at the Green function $G_{\infty}(x,y)$ on half space  $\mathbb{R}^n_+$:
$$ G_{\infty}(x,y) = \frac{A_{n,\alpha}}{s^{\frac{n-\alpha}{2}}}\left[1-
B\frac{1}{(t+s)^{\frac{n-2}{2}}}
\int_{0}^{\frac{s}{t}}\frac{(s-tb)^{\frac{n-2}{2}}}{b^{\alpha/2}(1+b)}db\right],$$
where $s =|x-y|^2$ and $t = 4x_ny_n$.
\bigskip

Next we establish the equivalence between problem (\ref{hpde}) and integral equation (\ref{int1}).

\begin{thm}\label{thm4.1}
Assume that $u$ is a locally bounded positive solution of
\begin{equation}
\left\{\begin{array}{ll}(-\Delta)^{\alpha/2}u(x)=u^{p}(x),&\mbox{in}\;\;\mathbb{R}^n_+,\\
u(x)=0, &\mbox{in}\;\;\;\; \mathbb{R}_-^n ;
\end{array}
\right.
\label{hpde1}
\end{equation}
and there exists a sequence $\{x^k\} \in \mathbb{R}_+^n$, such that
\begin{equation}
\frac{u(x^k)}{(x^k_n)^{\alpha/2}} \ra 0 .
\label{Ad4}
\end{equation}

Then it is also a solution of
\begin{equation}
u(x) = \int_{\mathbb{R}^n_+} G_{\infty} (x,y) u^p (y) dy ;
\label{inta}
\end{equation}
and vice versa.
\end{thm}

To prove the Theorem, we need the following Harnack inequality for $\alpha$-harmonic functions on domains with
boundaries, its consequences on half-spaces, and the uniqueness of $\alpha$-harmonic functions on half-spaces.

\begin{pro}(Boundary Harnack, see \cite{CaS} or \cite{Bo})\label{pro2}
Let $f$, $g:\;R^{n}\rightarrow R$ be two nonnegative functions such that $(-\Delta)^{s}f=(-\Delta)^{s}g=0$
in a domain $\Omega$. Suppose that $x_{0}\in\partial\Omega$, $f(x)=g(x)=0$ for
any $x\in B_{1}\backslash\Omega$, and $\partial\Omega\cap B_{1}$
is a Lipschitz graph in the direction of $x_{1}$ with Lipschitz constant less than $1$. Then
there is a constant $C$ depending only on dimension such that
\begin{equation}
\sup_{x\in \Omega\cap B_{\frac{1}{2}}}\frac{f(x)}{g(x)}\leq C\inf_{x\in \Omega\cap B_{\frac{1}{2}}}\frac{f(x)}{g(x)}
\end{equation}
for any $x,y\in B_{\frac{1}{2}}(x_{0})$.
\end{pro}

Based on this Harnack inequality, we derive the uniqueness of $\alpha$-harmonic functions on half spaces.

\begin{lem}\label{lem33}

Assume that $w$ is a nonnegative solution of
\begin{equation}
\left\{
  \begin{array}{ll}
   (-\Delta)^{\alpha/2}w=0,\;\;& x\in \mathbb{R}_{+}^{n}\\
w= 0, &x\in \mathbb{R}_{-}^{n}.
  \end{array}
\right.
\label{Ad6}
\end{equation}

Then there is a constant $c_o > 0$ such that for any two points $x=(x_1, \cdots, x_n)$ and $y=(y_1, \cdots, y_n)$ in  $\mathbb{R}_+^n$, we have
\begin{equation}
\frac{w(y)}{(y_n)^{\alpha/2}} \geq c_o \frac{w(x)}{(x_n)^{\alpha/2}}.
\label{Ad1}
\end{equation}

Consequently, we have either
\begin{equation}\label{f66}
w(x)\equiv 0,\;x\in \mathbb{R}^{n},
\end{equation}
or there exists a constant $a_o > 0$, such that
\begin{equation}
w(x) \geq a_o (x_n)^{\alpha/2}, \;\; \forall \; x\in \mathbb{R}_+^{n}.
\label{Ad2}
\end{equation}
\end{lem}

\emph{Proof.}

Let
$$
g(x) =\left\{
            \begin{array}{ll}
             (x_{n})^{\alpha/2} ,&x_{n}>0,\\
0,&x_{n}\leq0 .
            \end{array}
          \right.
$$
Then it is well known that $g(x)$ is a non-zero solution of problem (\ref{Ad6}).

We compare $w(x)$ with this $\alpha$-harmonic function $g(x)$. In Proposition \ref{pro2}, choose $\Omega = \mathbb{R}_+^n$.
It is easy to see that, by re-scaling, one can replace $B_{\frac{1}{2}}$ in Proposition \ref{pro2} by $B_{\frac{R}{2}}$
for any $R> 0$, and the constant $c$ is independent of $R$. Given any two points $x$ and $y$ in $\mathbb{R}_+^n$, choose
$R$ sufficiently large, such that $x, y \in B_{\frac{R}{2}}(0)$, then it follows from the Proposition that
\begin{equation}
\frac{w(y)}{(y_n)^{\alpha/2}} \geq c_o \frac{w(x)}{(x_n)^{\alpha/2}}.
\label{Ad3}
\end{equation}

Now, suppose $w$ is not identically zero. Then there exist a point $x^o \in \mathbb{R}_+^n$, such that $w(x^o) > 0$.
Hence by (\ref{Ad3}), we derive
$$ w(y) \geq c_o \frac{w(x^o)}{(x^o_n)^{\alpha/2}} (y_n)^{\alpha/2} \equiv a_o (y_n)^{\alpha/2} .$$

This completes the proof of the Lemma.

Now, we are ready to prove Theorem \ref{thm4.1}.
Assume that $u$ is a positive solution of (\ref{hpde1}). We first show that
\begin{equation}
\int_{\mathbb{R}_{+}^{n}}G_{\infty} (x,y) u^p (y) dy<\infty.
\end{equation}
Set
$$
v_{R}(x)=\int_{B_{R}(P_{R})}G_{R}(x,y)u^{p}(y)dy.
$$
From the local bounded-ness assumption on $u$, one can see that,  for each $R>0$, $v_R(x)$ is well-defined and is continuous. Moreover
$$
\left\{
  \begin{array}{ll}
   (-\Delta)^{\alpha/2}v_{R}(x)=u^{p}(x), & x\in B_{R}(P_{R}) \\
   v_{R}(x)=0,&x\not\in B_{R}(P_{R}).
  \end{array}
\right.
$$

Let $w_{R}(x)=u(x)-v_{R}(x)$, then
$$
\left\{
  \begin{array}{ll}
   (-\Delta)^{\alpha/2}w_{R}(x)=0, & x\in B_{R}(P_{R}), \\
   w_{R}(x)\geq0,&x\not\in B_{R}(P_{R}).
  \end{array}
\right.
$$

By the following Maximum Principle:  

\begin{pro} (Silvestre \cite{Si})
Let $\Omega \in \mathbb{R}^n$ be a bounded open set, and let $f$ be a lower-semicontinuous function in $\bar{\Omega}$ such that
$(- \lap)^{\alpha/2} f \geq 0$ in $\Omega$ and $f \geq 0$ in $\mathbb{R}^n \setminus \Omega.$ Then $f \geq 0$ in $\mathbb{R}^n$.
\label{proSi}
\end{pro}
We derive 
$$
w_{R}(x)\geq 0,\;\forall \;x\in B_{R}(P_{R}).
$$

Now letting  $R\rightarrow\infty$, we arrive at
$$
u(x)\geq \int_{\mathbb{R}_{+}^{n}}G_{\infty}(x,y)u^{p}(y)dy.
$$

Let
$$ v(x) = \int_{\mathbb{R}^n_+} G_{\infty} (x,y) u^p (y) dy .$$
Then
$$(- \lap )^{\alpha/2} v = u^p(x) , \;\; \forall \, x \in \mathbb{R}^n_+.$$

Set $w = u-v$, we have
\begin{equation}\label{90}
 \left\{\begin{array}{ll}
(-\lap )^{\alpha/2} w = 0 , \;w\geq 0,& x \in \mathbb{R}^n_+ ,\\
w (x) \equiv 0 , & x \in \mathbb{R}^n_- .
\end{array}
\right.
\end{equation}
Then by Lemma \ref{lem33}, we deduce
$$ w(x) \equiv 0 , \;\; \forall \, x \in \mathbb{R}^n .$$
Therefore,
$$ u(x) = v(x) = \int_{\mathbb{R}^n_+} G_{\infty} (x,y) u^p (y) dy .$$

On the other hand, assume that $u(x)$ is a solution of integral equation (\ref{inta}).
Applying $(-\lap)^{\alpha/2}$ to both sides, and on the right hand side, exchanging
it with the integral, we arrive at (\ref{hpde1}).

This completes the proof of the theorem.

\section{The Liouville Type Theorems in $\mathbb{R}^n_+$}
In this section, we prove the non-existence of positive solutions
under global and local integrability assumptions respectively and thus establish
Theorem \ref{mthm3} and \ref{mthm4}.

\subsection{Properties of the Green's Functions}
Let $\lambda$ be a positive real number and let the moving plane be
$$T_{\lambda}=\{x\in \mathbb{R}^n_+|\;x_{n}=\lambda\}.$$ We denote
$\Sigma_{\lambda}$ the region between the plane $x_{n}=0$ and the
plane $x_{n}=\lambda$. That is
$$
\Sigma_{\lambda}=\{x=(x_{1},\cdots,x_{n-1},x_{n})\in
\mathbb{R}^n_+|\;0<x_{n}<\lambda\}.
$$
Let
$$
x^{\lambda}=(x_{1},\cdots,x_{n-1},2\lambda-x_{n})
$$
be the reflection of the point
$x=(x_{1},\cdots,x_{n-1},x_{n})$ about the plane
$T_{\lambda}$, set
$$
\Sigma_{\lambda}^{C}=\mathbb{R}^n_+\backslash\Sigma_{\lambda}
$$
the complement of $\Sigma_{\lambda} $, and write
$$
u_{\lambda}(x)=u(x^{\lambda})\;\;\mbox{and}\;\;w_{\lambda}(x)=u_{\lambda}(x)-u(x).
$$

First we derive the properties of $G_{\infty}(x,y)$.

\begin{lem}\label{alemm1}
(i) For any $x,y\in \Sigma_{\lambda}$, $x\neq y$, we have
\begin{equation}\label{a6}
G_{\infty}(x^{\lambda},y^{\lambda})>\max\{G_{\infty}(x^{\lambda},y),
G_{\infty}(x,y^{\lambda})\}
\end{equation}
and
\begin{equation}\label{a7}
G_{\infty}(x^{\lambda},y^{\lambda})-G_{\infty}(x,y)>|G_{\infty}(x^{\lambda},y)-G_{\infty}(x,y^{\lambda})|.
\end{equation}
(ii) For any $x\in \Sigma_{\lambda}$, $y\in \Sigma_{\lambda}^{C}$, it holds
$$
G_{\infty}(x^{\lambda},y)> G_{\infty}(x,y).
$$
\end{lem}
\textbf{Proof.}
$(i)$   Set $\varphi(x,y)=4x_{n}y_{n}$. Obviously, we have
\begin{eqnarray}\label{a010}
\varphi(x^{\lambda},y^{\lambda})&>&\max\{\varphi(x,y^{\lambda}),\varphi(x^{\lambda},y)\}\\
&\geq&\min\{\varphi(x,y^{\lambda}),\varphi(x^{\lambda},y)\}\\
\label{a0.13}
&>&\varphi(x,y).
\end{eqnarray}
Let $t=\varphi(x,y)>0$, $s=d(x,y)=|x-y|^{2}>0$,
$$
I_{\infty}(s,t)=\int_{0}^{\frac{s}
{t}}\frac{(\frac{1-\frac{t}{s}b}{1+\frac{t}{s}})
^{\frac{n-2}{2}}}{b^{\alpha/2}(1+b)}db
$$
and
$$
H_{\infty}(s,t)=\frac{1}{s^{\frac{n-\alpha}{2}}}\left[1-B I_{\infty}(s,t)\right].
$$
Then similar to the $\textbf{B}_{1}$ case, we have
\begin{equation}\label{a13}
\frac{\partial H_{\infty}}{\partial s}=-\frac{n-\alpha}{2s^{\frac{n-\alpha}{2}+1}}\left[1-BI_{\infty}(s,t)\right]+\frac{-B}{s^{\frac{n-\alpha}{2}}}\frac{\partial I_{\infty}}{\partial s}<0.
\end{equation}
\begin{equation}\label{a14}
\frac{\partial H_{\infty}}{\partial t}=\frac{-B}{s^{\frac{n-\alpha}{2}}}\frac{\partial I_{\infty}}{\partial t}>0.
\end{equation}
\begin{equation}\label{a81}
\frac{\partial^{2}H_{\infty}}{\partial t\partial s}<0.
\end{equation}
Since
\begin{equation}\nonumber
G_{\infty}(x^{\lambda},y^{\lambda})=A_{n,\alpha}H_{\infty}(d(x^{\lambda},y^{\lambda}),\varphi
(x^{\lambda},y^{\lambda})),
\end{equation}
\begin{equation}\nonumber
G_{\infty}(x^{\lambda},y)=A_{n,\alpha}H_{\infty}(d(x^{\lambda},y),\varphi
(x^{\lambda},y)),
\end{equation}
\begin{equation}\nonumber
G_{\infty}(x,y^{\lambda})=A_{n,\alpha}H_{\infty}(d(x,y^{\lambda}),\varphi
(x,y^{\lambda})),
\end{equation}
\begin{equation}\nonumber
G_{\infty}(x,y)=A_{n,\alpha}H_{\infty}(d(x,y),\varphi (x,y)),
\end{equation}
then (\ref{a6}) is a direct consequence of (\ref{a010}), (\ref{a13}), and  (\ref{a14}).

Next, we prove (\ref{a7}). From
(\ref{a010})-(\ref{a0.13}) and (\ref{a81}), we deduce
\begin{eqnarray}\nonumber
&&G_{\infty}(x^{\lambda},y^{\lambda})-G_{\infty}(x,y)\\
\nonumber
&=&A_{n,\alpha}\int_{\varphi(x,y)}^{\varphi(x^{\lambda},y^{\lambda})}\frac{\partial
H_{\infty}(d(x,y),t)}{\partial t}dt\\
\nonumber
&>&A_{n,\alpha}\int_{\varphi(x,y)}^{\varphi(x^{\lambda},y^{\lambda})}\frac{\partial
H_{\infty}(d(x^{\lambda},y),t)}{\partial t}dt\\
\nonumber &\geq&
A_{n,\alpha}\int_{\varphi(x,y^{\lambda})}^{\varphi(x^{\lambda},y)}\frac{\partial
H_{\infty}(d(x^{\lambda},y),t)}{\partial t}dt\\
\nonumber
&=&A_{n,\alpha}[H_{\infty}(d(x^{\lambda},y),\varphi(x^{\lambda},y_{n}))-H_{\infty}(d(x,y^{\lambda}),\varphi(x,y^{\lambda}))]\\
&=& G_{\infty}(x^{\lambda},y)-G_{\infty}(x,y^{\lambda}).
\end{eqnarray}

Similarly, one can show that
$$G_{\infty} (x^{\lambda}, y^{\lambda}) - G_{\infty} (x,y) > G_{\infty}(x, y^{\lambda}) - G_{\infty}(x^{\lambda},y) .$$
These imply (\ref{a7}).

$(ii)$ For $x\in \Sigma_{\lambda}$ and $y\in \Sigma_{\lambda}^{C}$, we have
\begin{equation}\label{a52}
d(x^{\lambda},y)<d(x,y),\;\;\varphi(x^{\lambda},y)>\varphi(x,y).
\end{equation}
By (\ref{a13})-(\ref{a14}), and (\ref{a52}), we get
$$
G_{\infty}(x^{\lambda},y)> G_{\infty}(x,y).
$$
This completes the proof of Lemma \ref{alemm1}.

The following lemma is a key ingredient in our integral estimates and the proof is similar to that of Lemma \ref{lemm2}, which we  omit here.
\begin{lem}\label{lemm52}
For any $x\in\Sigma_{\lambda}$, it holds
$$
u(x)-u_{\lambda}(x)\leq \int_{\Sigma_{\lambda}}\left[
G_{\infty}(x^{\lambda},y^{\lambda})-G_{\infty}(x,y^{\lambda})\right]\left[u^{p}(y)-u^{p}_{\lambda}(y)\right]dy.
$$
\end{lem}

\subsection{Non-existence under Global Integrability Assumption}

In this subsection, we prove
\begin{thm}\label{bthm6}
  Assume $p>\frac{n}{n-\alpha}$.
If $u\in L^{\frac{n(p-1)}{\alpha}}(\mathbb{R}^n_+)$ is a non-negative solution of
\begin{equation}\label{5a1}
u(x)=\int_{\mathbb{R}^n_+}G_{\infty}(x,y)u^{p}(y)dy,
\end{equation}
then $u(x)\equiv 0$.
\end{thm}

From the following result of Silvestre \cite{Si}, one can see that a nonnegative solution $u$ is either strictly positive or
identically zero in $\mathbb{R}^n$.

\begin{pro} Let $\Omega \Subset \mathbb{R}^n$ be an open set, and let $u$ be a lower-semi-continuous
function in $\bar{\Omega}$ such that
$$ (-\lap)^{\alpha/2} u \geq 0 \; \mbox{ and } \; u \geq 0 \; \mbox{ in } \; \mathbb{R}^n \setminus \Omega . $$
Then $u \geq 0$ in $\mathbb{R}^n$. Moreover, if $u(x)=0$ for some point inside $\Omega$, then $u \equiv 0$ in all
$\mathbb{R}^n$.
\label{pro5.1}
\end{pro}

By virtue of this proposition, without loss of generality, we may assume that $u>0$ in $\mathbb{R}^n_+$ and derive a contradiction.

We divide the proof into two steps. In the first step, we start from
the very low end of our region $\mathbb{R}^n_+$, i.e. near $x_{n}=0$.
We will show that for $\lambda$ sufficiently small,
\begin{equation}\label{k21}
w_{\lambda}(x)=u_{\lambda}(x)-u(x)\geq 0,\;\;a.e.\;\forall
x\in\Sigma_{\lambda}.
\end{equation}
In the second step, we will move our plane $T_{\lambda}$ up in the
positive $x_{n}$ direction as long as the inequality
(\ref{k21}) holds to show that $u(x)$ is monotone increasing in $x_n$
and thus derive a contradiction.

\emph{Step 1.} Define
$$
\Sigma_{\lambda}^{-}=\{x\in\Sigma_{\lambda}| \;\;w_{\lambda}(x)<0
\}.
$$
We show that for $\lambda$ sufficiently small,
$\Sigma_{\lambda}^{-}$ must be measure zero. In fact, for any $x\in \Sigma _{\lambda}^{-}$, by the Mean Value Theorem and Lemma \ref{lemm52}, we obtain
\begin{eqnarray}
\nonumber
0&<&u(x)-u_{\lambda}(x)\\
\nonumber
 &\leq &\int_{\Sigma_{\lambda}}\left[
G_{\infty}(x^{\lambda},y^{\lambda})-
G_{\infty}(x,y^{\lambda})\right]\left[u^{p}(y)-u_{\lambda}^{p}(y)\right]dy\\
\nonumber
 &=&\int_{\Sigma^{-}_{\lambda}}\left[
G_{\infty}(x^{\lambda},y^{\lambda})-
G_{\infty}(x,y^{\lambda})\right]\left[u^{p}(y)-u_{\lambda}^{p}(y)\right]dy\\
\nonumber
&&+\int_{\Sigma_{\lambda}\backslash\Sigma^{-}_{\lambda}}\left[
G_{\infty}(x^{\lambda},y^{\lambda})-
G_{\infty}(x,y^{\lambda})\right]\left[u^{p}(y)-u_{\lambda}^{p}(y)\right]dy\\
\nonumber
 &\leq&\int_{\Sigma^{-}_{\lambda}}\left[
G_{\infty}(x^{\lambda},y^{\lambda})-
G_{\infty}(x,y^{\lambda})\right]\left[u^{p}(y)-u_{\lambda}^{p}(y)\right]dy\\
\nonumber
 &\leq&\int_{\Sigma^{-}_{\lambda}}
G_{\infty}(x^{\lambda},y^{\lambda})\left[u^{p}(y)-u_{\lambda}^{p}(y)\right]dy\\
\nonumber
 &=&p\int_{\Sigma^{-}_{\lambda}}
G_{\infty}(x^{\lambda},y^{\lambda})\psi^{p-1}_{\lambda}(y)[u(y)-u_{\lambda}(y)]dy\\
\label{k1}
&\leq& p\int_{\Sigma^{-}_{\lambda}}
G_{\infty}(x^{\lambda},y^{\lambda})u^{p-1}(y)[u(y)-u_{\lambda}(y)]dy,
\end{eqnarray}
where $\psi_{\lambda}(y)$ is valued between $u(y)$ and
$u_{\lambda}(y)$ and hence on $\Sigma_{\lambda}^{-}$, we have
$$
0\leq u_{\lambda}(y)\leq\psi_{\lambda}(y)\leq u(y).
$$
By the expression of $G_{\infty}(x,y)$, it is easy to see
$$
G_{\infty}(x,y)
\leq \frac{A_{n,\alpha}}{|x-y|^{n-\alpha}}.
$$
 Now (\ref{k1}) implies
\begin{equation}\label{k2}
0<u(x)-u_{\lambda}(x)\leq\int_{\Sigma^{-}_{\lambda}}
\frac{C}{|x-y|^{n-\alpha}}|u^{p-1}(y)||u(y)-u_{\lambda}(y)|dy.
\end{equation}
We apply the Hardy-Littlewood-Sobolev inequality (\ref{f1}) and H\"{o}lder
inequality to (\ref{k2}) to obtain, for any $q>\frac{n}{n-\alpha}$,
\begin{eqnarray}
\nonumber
\|w_{\lambda}\|_{L^{q}(\Sigma_{\lambda}^{-})}&\leq&
C\|u^{p-1}w_{\lambda}\|_{L^{\frac{nq}{n+\alpha
q}}(\Sigma_{\lambda}^{-})}\\
\label{k3}
 &\leq& C\|u^{p-1}\|_{L^{\frac{n}{\alpha}}(\Sigma_{\lambda}^{-})}
\|w_{\lambda}\|_{L^{q}(\Sigma_{\lambda}^{-})}.
\end{eqnarray}
Note here we can choose $q=\frac{n(p-1)}{\alpha}$, then by our assumption $p > \frac{n}{n-\alpha}$, we have
$q > \frac{n}{n-\alpha}$ and $w_{\lambda} \in L^q(\mathbb{R}^n)$.

Since $u\in L^{\frac{n(p-1)}{\alpha}}(\mathbb{R}^n_+)$, we can choose
sufficiently small positive $\lambda$ such that
\begin{equation}\label{k4}
C\|u^{p-1}\|_{L^{\frac{n}{\alpha}}(\Sigma_{\lambda}^{-})}=C\{\int_{\Sigma^{-}_{\lambda}}
u^{\frac{n(p-1)}{\alpha}}(y)dy\}^{\frac{\alpha}{n}}\leq\frac{1}{2}.
\end{equation}
By (\ref{k3})-(\ref{k4}), we have
$$
\|w_{\lambda}\|_{L^{q}(\Sigma_{\lambda}^{-})}=0,
$$
and therefore $\Sigma_{\lambda}^{-}$ must be measure zero.

\emph{Step 2.}  Now we start from
such small $\lambda$ and move the plane $T_{\lambda}$ up as
long as (\ref{k21}) holds.

Define
$$
\lambda_{0}=\sup\{\lambda|\;\;w_{\rho}(x)\geq 0,\; \rho\leq\lambda,\;
\forall x\in\Sigma_{\rho}\}.
$$
We will prove
\begin{equation}\label{d6a}
\lambda_{0}=+\infty.
\end{equation}

Suppose in the contrary that $\lambda_{0}<\infty$, we will show that $u(x)$ is symmetric about the plane
$T_{\lambda_{0}}$, i.e.
\begin{equation}\label{a3a}
w_{\lambda_{0}}\equiv 0,\;\;a.e.\;\;\forall x\in
\Sigma_{\lambda_{0}}.
\end{equation}
This will contradict the strict positivity of $u$.

Suppose (\ref{a3a}) does not hold. Then for such a $\lambda_{0}$, we have $w_{\lambda_{0}}\geq 0$, but $w_{\lambda_{0}}\not\equiv0$
$a.e.$ on $\Sigma_{\lambda_{0}}$. We show that the plane can be moved
further up.  More precisely, there exists an
$\epsilon>0$ such that for all $\lambda \in
[\lambda_{0},\lambda_{0}+\epsilon)$
\begin{equation}
w_{\lambda} (x) \geq 0 , \;\;a.e.\;\;\mbox{on}\;\;\Sigma_{\lambda}.
\label{5.2c}
\end{equation}
To verify this, we will again resort to inequality (\ref{k3}). If one can show that for
$\epsilon$ sufficiently small so that for all $\lambda$ in
$[\lambda_{0},\lambda_{0}+\epsilon),$ holds
\begin{equation}\label{a30a}
C\left\{\int_{\Sigma^{-}_{\lambda}}
u^{\frac{n(p-1)}{\alpha}}(y)dy\right\}^{\frac{\alpha}{n}}\leq\frac{1}{2},
\end{equation}
then by
(\ref{k3}) and (\ref{a30a}), we have
$\|w_{\lambda}\|_{L^{q}(\Sigma_{\lambda}^{-})}=0$, and
therefore $\Sigma_{\lambda}^{-}$ must be measure zero. Hence, for
these values of $\lambda>\lambda_{0}$, we have (\ref{5.2c}).
This contradicts the definition of $\lambda_{0}$. Therefore
(\ref{a3a}) must hold.

We postpone the proof of (\ref{a30a}) for a moment.

By (\ref{a3a}), we derive that the plane $x_{n}=2\lambda_{0}$ is the symmetric image of
the boundary $\partial \mathbb{R}^n_+$ with respect to the plane $T_{\lambda_{0}}$, and hence $u(x)=0$ when $x$
is on the plane $x_{n}=2\lambda_{0}$. This contradicts our assumption $u(x)>0$ in $\mathbb{R}^n_+$. Therefore, (\ref{d6a}) must be valid.

We have proved that the positive solution of (\ref{5a1}) is monotone increasing with respect to $x_{n}$, and this contradicts $u\in L^{\frac{n(p-1)}{\alpha}}(\mathbb{R}^n_+)$. Hence the positive solutions of (2) do not exist.

Now we verify inequality (\ref{a30a}). For any small $\eta>0$, we
can choose $R$ sufficiently large so that
\begin{equation}\label{a97a}
\left(\int_{\mathbb{R}^n_+\backslash
B_{R}}u^{\frac{n(p-1)}{\alpha}}(y)dy\right)^{\frac{\alpha}{n}}<\eta.
\end{equation}
We fix this $R$ and then show that the measure of
$\Sigma_{\lambda}^{-}\cap B_{R}$ is sufficiently small for
$\lambda$ close to $\lambda_{0}$.  First, we have
\begin{equation}\label{a25a}
w_{\lambda_{0}}(x)>0
\end{equation}
in the interior of $\Sigma_{\lambda_{0}}$.

Indeed, since
\begin{eqnarray}
\nonumber
u_{\lambda}(x)-u(x)&\geq&\int_{\Sigma_{\lambda}}\left[G_{\infty}(x^{\lambda},y^{\lambda})
-G_{\infty}(x,y^{\lambda})\right]\left[u^{p}_{\lambda}(y))-u^{p}(y)\right]dy\\
\nonumber
&\;\;\;\;+&\int_{\Sigma_{\lambda}^{C}\backslash\widetilde{\Sigma}_{\lambda}}
\left[G_{\infty}(x^{\lambda},y)-G_{\infty}(x,y)\right]u^{p}(y)dy\\
\label{e8a}
&\geq& \int_{\Sigma_{\lambda}^{C}\backslash\widetilde{\Sigma}_{\lambda}}
\left[G_{\infty}(x^{\lambda},y)-G_{\infty}(x,y)\right]u^{p}(y)dy.
\end{eqnarray}
If (\ref{a25a}) is violated, there exists some point
$x_{0}\in\Sigma_{\lambda_{0}}$ such that
$u(x_{0})=u_{_{\lambda_{0}}}(x_{0})$. And then by (\ref{e8a}) and Lemma \ref{alemm1} (ii), we obtain
\begin{equation}\label{a26}
u(y)\equiv 0,\;\;\forall y\in\Sigma_{\lambda_{0}}^{C}\backslash\widetilde{\Sigma}_{\lambda_{0}} .
\end{equation}
This is a contradiction with our assumption that $u>0$. Therefore
(\ref{a25a}) holds.

For any $\gamma>0$, let
\begin{equation}\label{44a}
E_{\gamma}=\{x\in \Sigma_{\lambda_{0}}\cap
B_{R}|\;w_{\lambda_{0}}(x)>\gamma\},\;\;F_{\gamma}=(\Sigma_{\lambda_{0}}\cap
B_{R})\backslash E_{\gamma}.
\end{equation}
It is obviously that
$$
\lim_{\gamma\rightarrow 0}\mu(F_{\gamma})=0.
$$
For $\lambda>\lambda_{0}$, let
$$
D_{\lambda}=(\Sigma_{\lambda}\backslash\Sigma_{\lambda_{0}})\cap
B_{R}.
$$
Then it is easy to see that
\begin{equation}
\label{42a}
(\Sigma^{-}_{\lambda}\cap B_{R})\subset (\Sigma^{-}_{\lambda}\cap
E_{\gamma})\cup F_{\gamma}\cup D_{\lambda}.
\end{equation}
Apparently, the measure of $D_{\lambda}$ is small for $\lambda$
close to $\lambda_{0}$. We show that the measure of
$\Sigma^{-}_{\lambda}\cap E_{\gamma}$ can also be sufficiently small as
$\lambda$ close to $\lambda_{0}$. In fact, for any
$x\in\Sigma^{-}_{\lambda}\cap E_{\gamma}$, we have
$$
w_{\lambda}(x)=u_{\lambda}(x)-u(x)=u_{\lambda}(x)-u_{\lambda_{0}}(x)+u_{\lambda_{0}}(x)-u(x)<0.
$$
Hence
$$
u_{\lambda_{0}}(x)-u_{\lambda}(x)>w_{\lambda_{0}}(x)>\gamma.
$$
It follows that
\begin{equation}\label{41a}
(\Sigma^{-}_{\lambda}\cap E_{\gamma})\subset G_{\gamma}\equiv \{x\in
B_{R}|\;\;u_{\lambda_{0}}(x)-u_{\lambda}(x)>\gamma\}.
\end{equation}
By the well-known Chebyshev inequality, we have
\begin{eqnarray}\nonumber
\mu(G_{\gamma})&\leq&\frac{1}{\gamma^{p+1}}\int_{G_{\gamma}}|u_{\lambda_{0}}(x)-u_{\lambda}(x)|^{p+1}dx\\
&\leq&
\frac{1}{\gamma^{p+1}}\int_{B_{R}}|u_{\lambda_{0}}(x)-u_{\lambda}(x)|^{p+1}dx.
\end{eqnarray}
For each fixed $\gamma$, as $\lambda$ close to $\lambda_{0}$, the
right hand side of the above inequality can be made as small as we
wish. Therefore by (\ref{42a}) and (\ref{41a}), the measure of
$\Sigma^{-}_{\lambda}\cap B_{R}$ can also be made sufficiently
small. Combining this with (\ref{a97a}), we obtain (\ref{a30a}).
This completes the proof of
Theorem \ref{bthm6}.

\subsection{Non-existence under Weaker Conditions}

In this section, we will use a proper Kelvin type transforms and derive
non-existence of positive solutions in $\mathbb{R}^n_+$ under much weaker conditions, i.e. the solution $u$ is only locally bounded or, in the critical case, 
only locally integrable.

Because there is no explicit global integrability assumptions on the solution $u$, we cannot directly carry on the method
of moving planes on $u$. To overcome this difficulty, we employ Kelvin type transforms.

For $z^{0}\in\partial \mathbb{R}^n_+$, let
\begin{equation}\label{h62}
\bar{u}_{z^0}(x)=\frac{1}{|x-z^{0}|^{n-\alpha}}u\left(\frac{x-z^{0}}{|x-z^{0}|^{2}}+z^{0}\right)
\end{equation}
be the Kelvin type transform of $u$ centered at $z^0$.

Through a straight forward calculation, we have
\begin{eqnarray}
\nonumber
\bar{u}_{z^0}(x)
&=&\frac{1}{|x-z^{0}|^{n-\alpha}}\int_{\mathbb{R}^n_+}G_{\infty}\left(\frac{x-z^{0}}{|x-z^{0}|^{2}}+z^{0},y \right)u^{p}(y)dy \nonumber \\
\label{53}
&=&\int_{\mathbb{R}^n_+}G_{\infty}(x,y)\frac{\bar{u}_{z^0}^{p}(y)}{|y-z^{0}|^{\beta}}dy,\;\;\forall\;
x\in \mathbb{R}^n_+\backslash B_{\epsilon}(z^{0}),\;\epsilon>0,
\end{eqnarray}
where $\frac{n}{n-\alpha}<p\leq\tau$, $\beta=(n-\alpha)(\tau-p)\geq 0$, and $\tau=\frac{n+\alpha}{n-\alpha}$.

We consider critical case and subcritical case separately.

\textbf{(i)} {\em The critical case} $p=\tau=\frac{n+\alpha}{n-\alpha}$. If $u(x)$ is a solution of
\begin{equation}\label{66}
u(x)=\int_{\mathbb{R}^n_+}G_{\infty}(x,y)u^{\tau}(y)dy,
\end{equation}
then $\bar{u}_{z^0}$ is also a solution of (\ref{66}). Since $u\in L_{loc}^{\frac{2n}{n-\alpha}}$, for any domain $\Omega$ that is a positive distance away from $z^{0}$, we have
\begin{equation}\label{b85}
\int_{\Omega}\bar{u}_{z^0}^{\frac{2n}{n-\alpha}}(y)dy<\infty.
\end{equation}

We consider two possibilities.

\textbf{Possibility  1.} If there is a $z^{0}=(z^{0}_{1},\cdots,z^{0}_{n-1},0)\in \partial \mathbb{R}^n_+$
such that $\bar{u}_{z^0}(x)$ is bounded near $z^{0}$, then by (\ref{h62}), we obtain
$$
u(y)=\frac{1}{|y-z^{0}|^{n-\alpha}}\bar{u}_{z^0}\left(\frac{y-z^{0}}{|y-z^{0}|^{2}}+z^{0} \right).
$$
And we further deduce
\begin{equation}\label{63}
u(y)=O \left(\frac{1}{|y|^{n-\alpha}}\right), \;\;\mbox{as}\;|y|\rightarrow\infty.
\end{equation}
Since $u\in L_{loc}^{\frac{2n}{n-\alpha}}(\mathbb{R}^n_+)$, together with (\ref{63}), we have
\begin{equation}\label{a85}
\int_{\mathbb{R}^n_+}
u^{\frac{2n}{n-\alpha}}(y)dy<\infty.
\end{equation}

In this situation, we still carry on the moving of planes on $u$. Going through exactly the same arguments as in the proof of Theorem \ref{bthm6}, we derive the non-existence of positive solutions for (\ref{inta}).

\textbf{Possibility 2.}  For all $z^{0}=(z^{0}_{1},\cdots,z^{0}_{n-1},0)\in \partial \mathbb{R}^n_+$,
 $\bar{u}_{z^0}(x)$ are unbounded near $z^{0}$. Then for each $z^0$,  we will carry on the moving planes on $\bar{u}_{z^0}$ in $\mathbb{R}^{n-1}$ to prove that it is rotationally symmetric about
the line passing through $z^{0}$ and parallel to the $x_{n}$-axis. From this, we will deduce that $u$ is independent of the first $n-1$ variables $x_1, \cdots, x_{n-1}$. That is, $u = u(x_n)$, which as we will show, contradicts the finiteness of the integral
$$ \int_{\mathbb{R}^n_+} G_{\infty}(x,y) u^p (y) d y .$$

In this situation, since we only need to deal with $\bar{u}_{z^0}$, for simplicity, we denote it by $\bar{u}$.
For a given real number $\lambda$, define
$$
\hat{\Sigma}_{\lambda}=\{x=(x_{1},\cdots,x_{n})\in
\mathbb{R}^n_+ \mid x_{1}<\lambda\}
$$
and let
$$
x^{\lambda}=(2\lambda-x_{1},x_{2},\cdots,x_{n}).
$$

For $x,y\in \hat{\Sigma}_{\lambda}$, $x\neq y$, by (\ref{a13}), it is easy to see
\begin{eqnarray}
G_{\infty}(x,y)=G_{\infty}(x^{\lambda},y^{\lambda}),\;\;G_{\infty}(x^{\lambda},y)&=&G_{\infty}(x,y^{\lambda}), \nonumber
 \\ 
 \mbox{and}\;\;G_{\infty}(x^{\lambda},y^{\lambda})&>&G_{\infty}(x,y^{\lambda}).
 \label{54}
\end{eqnarray}

Obviously, we have
$$
\bar{u}(x)=\int_{\hat{\Sigma}_{\lambda}}G_{\infty}(x,y)\bar{u}^{\tau}(y)dy+\int_{\hat{\Sigma}_{\lambda}}G_{\infty}
(x,y^{\lambda})\bar{u}^{\tau}_{\lambda}(y)dy
$$
$$
\bar{u}(x^{\lambda})=\int_{\hat{\Sigma}_{\lambda}}G_{\infty}(x^{\lambda},y)\bar{u}^{\tau}(y)dy
+\int_{\hat{\Sigma}_{\lambda}}G_{\infty}(x^{\lambda},y^{\lambda})\bar{u}^{\tau}_{\lambda}(y)dy.
$$

By (\ref{54}), one derives
\begin{eqnarray}
\nonumber
&&\bar{u}(x)-\bar{u}(x^{\lambda})=\int_{\hat{\Sigma}_{\lambda}}\left[G_{\infty}(x,y)-G_{\infty}(x^{\lambda},y)\right]\bar{u}^{\tau}(y)dy\\
\nonumber
&\;\;\;\;
+&\int_{\hat{\Sigma}_{\lambda}}
\left[G_{\infty}(x,y^{\lambda})-G_{\infty}(x^{\lambda},y^{\lambda})\right]\bar{u}^{\tau}_{\lambda}(y)dy\\
\label{14}
&=&\int_{\hat{\Sigma}_{\lambda}}\left[G_{\infty}(x,y)-G_{\infty}(x^{\lambda},y)\right]\left[\bar{u}^{\tau}(y)-\bar{u}^{\tau}_{\lambda}(y)\right]dy.
\end{eqnarray}

In the first step,
we will show that for $\lambda$ sufficiently negative,
\begin{equation}\label{b21}
w_{\lambda}(x) \equiv \bar{u}_{\lambda}(x)-\bar{u}(x)\geq 0,\;\;a.e.\;\forall
x\in\hat{\Sigma}_{\lambda}.
\end{equation}
In the second step, we deduce that $\hat{T}$ can be move to the right all the way to $z^{0}$.  And furthermore, we derive $w_{z_{1}^{0}}\equiv 0,$ $\forall x\in \hat{\Sigma}_{z_{1}^{0}}$.

\emph{Step 1.}  Define
$$
\hat{\Sigma}_{\lambda}^{-}=\{x\in\hat{\Sigma}_{\lambda}\backslash B_{\epsilon}((z^{0})^{\lambda})|\;w_{\lambda}(x)<0
\},
$$
where $(z^{0})^{\lambda}$ is the reflection of $z^{0}$ about the plane $\hat{T}_{\lambda}=\{x\in \mathbb{R}^n_+
|\;x_{1}=\lambda\}$.
We show that for $\lambda$ sufficiently negative,
$\hat{\Sigma}_{\lambda}^{-}$ must be measure zero. In fact, by (\ref{14}), (\ref{54}), and the Mean Value Theorem, we obtain, for $x\in\hat{\Sigma}_{\lambda}^{-}$,
\begin{eqnarray}
\nonumber
0&<&\bar{u}(x)-\bar{u}_{\lambda}(x)\\
\nonumber
&=&\int_{\hat{\Sigma}^{-}_{\lambda}}[
G_{\infty}(x,y)-G_{\infty}(x,y^{\lambda})][\bar{u}^{\tau}(y)-\bar{u}_{\lambda}^{\tau}(y)]dy\\
\nonumber
&\;\;+&\int_{\hat{\Sigma}_{\lambda}\backslash\hat{\Sigma}^{-}_{\lambda}}[G_{\infty}(x,y)
-G_{\infty}(x,y^{\lambda})][\bar{u}^{\tau}(y)-\bar{u}_{\lambda}^{\tau}(y)]dy\\
\nonumber
 &\leq&\int_{\hat{\Sigma}^{-}_{\lambda}}[
G_{\infty}(x,y)-G_{\infty}(x,y^{\lambda})][\bar{u}^{\tau}(y)-\bar{u}_{\lambda}^{\tau}(y)]dy\\
\nonumber
 &\leq&\int_{\hat{\Sigma}^{-}_{\lambda}}
G_{\infty}(x,y)[\bar{u}^{\tau}(y)-\bar{u}_{\lambda}^{\tau}(y)]dy\\
\nonumber
 &=&\tau\int_{\hat{\Sigma}^{-}_{\lambda}}
G_{\infty}(x,y)\psi^{\tau-1}_{\lambda}(y)[\bar{u}(y)-\bar{u}_{\lambda}(y)]dy\\
\nonumber
&\leq& \tau\int_{\hat{\Sigma}^{-}_{\lambda}}
G_{\infty}(x,y)\bar{u}^{\tau-1}(y)[\bar{u}(y)-\bar{u}_{\lambda}(y)]dy\\
\label{b22}
&\leq&\int_{\hat{\Sigma}^{-}_{\lambda}}
\frac{C}{|x-y|^{n-\alpha}}|\bar{u}^{\tau-1}(y)||\bar{u}(y)-\bar{u}_{\lambda}(y)|dy.
\end{eqnarray}
We apply the Hardy-Littlewood-Sobolev inequality (\ref{f1}) and H\"{o}lder
inequality to (\ref{b22}) to obtain, for any $q>\frac{n}{n-\alpha}$,
\begin{eqnarray}
\nonumber
\|w_{\lambda}\|_{L^{q}(\hat{\Sigma}_{\lambda}^{-})}&\leq&
C\|\bar{u}^{\tau-1}w_{\lambda}\|_{L^{\frac{nq}{n+\alpha
q}}(\hat{\Sigma}_{\lambda}^{-})}\\
\label{b24}
 &\leq &C\|\bar{u}^{\tau-1}\|_{L^{\frac{n}{\alpha}}(\hat{\Sigma}_{\lambda}^{-})}
\|w_{\lambda}\|_{L^{q}(\hat{\Sigma}_{\lambda}^{-})}.
\end{eqnarray}
By (\ref{b85}), we can choose $N$
sufficiently large such that for $\lambda\leq-N$,
\begin{equation}\label{a10}
C\|\bar{u}^{\tau-1}\|_{L^{\frac{n}{\alpha}}(\hat{\Sigma}_{\lambda}^{-})}=C\left\{\int_{\hat{\Sigma}^{-}_{\lambda}}
\bar{u}^{\frac{2n}{n-\alpha}}(y)dy\right\}^{\frac{\alpha}{n}}\leq\frac{1}{2}.
\end{equation}
Now inequalities (\ref{b24}) and (\ref{a10}) imply
$$
\|w_{\lambda}\|_{L^{q}(\hat{\Sigma}_{\lambda}^{-})}=0,
$$
and therefore $\hat{\Sigma}_{\lambda}^{-}$ must be measure zero.

\emph{Step 2.}  ({ \em Move the plane to the limiting position to derive\emph{}
symmetry.} )

Inequality (\ref{b21}) provides a starting point to move the plane
$\hat{T}_{\lambda}$. Now we start from
the neighborhood of $x_{1}=-\infty$ and move the plane to the right as
long as (\ref{b21}) holds to the limiting position. Define
$$
\lambda_{0}=\sup\{\lambda\leq z_{1}^{0}|\; w_{\rho}(x)\geq 0,\; \rho\leq\lambda,\;
\forall x\in\hat{\Sigma}_{\rho}\}.
$$
We prove that $\lambda_{0}\geq z_{1}^{0} - \epsilon$. On the contrary, suppose that $\lambda_{0}<z_{1}^{0} -\epsilon$.
We will show that $\bar{u}(x)$ is symmetric about the plane
$T_{\lambda_{0}}$, i.e.
\begin{equation}\label{b3}
w_{\lambda_{0}}\equiv 0,\;\;a.e.\;\;\forall x\in
\hat{\Sigma}_{\lambda_{0}}\backslash B_{\epsilon}((z^{0})^{\lambda_{0}}).
\end{equation}
Suppose (\ref{b3}) is not true, then  for such  $\lambda_{0} <z_{1}^{0} -\epsilon$, we have
$$w_{\lambda_{0}}\geq 0 , \mbox{ but } w_{\lambda_{0}}\not\equiv0
\mbox{ a.e. on } \hat{\Sigma}_{\lambda_{0}}\backslash B_{\epsilon}((z^{0})^{\lambda_{0}}).$$ We show that the plane can be
moved further to the right. More precisely, there exists a
$\zeta>0$ such that for all $\lambda\in[\lambda_{0},\lambda_{0}+\zeta)$
$$w_{\lambda}(x) \geq 0 , \;\;a.e\;\;\mbox{on}\;\;\hat{\Sigma}_{\lambda}\backslash B_{\epsilon}((z^{0})^{\lambda}).
$$
This will contradicts the definition of $\lambda_0$.

By inequality (\ref{b24}), we have
\begin{equation}\label{b29}
\|w_{\lambda}\|_{L^{q}(\hat{\Sigma}_{\lambda}^{-})}\leq
C\left\{\int_{\hat{\Sigma}^{-}_{\lambda}}
\bar{u}^{\frac{2n}{n-\alpha}}(y)dy\right\}^{\frac{\alpha}{n}}
\|w_{\lambda}\|_{L^{q}(\hat{\Sigma}_{\lambda}^{-})}.
\end{equation}
Similar to the proof of (\ref{a30a}), we can
choose $\zeta$ sufficiently small so that for all $\lambda\in[\lambda_{0},\lambda_{0}+\zeta)$,
\begin{equation}\label{b30}
C\{\int_{\hat{\Sigma}^{-}_{\lambda}}
\bar{u}^{\frac{2n}{n-\alpha}}(y)dy\}^{\frac{\alpha}{n}}\leq\frac{1}{2}.
\end{equation}
We postpone the proof of this inequality for a moment. Now by
(\ref{b29}) and (\ref{b30}), we have
$\|w_{\lambda}\|_{L^{q}(\hat{\Sigma}_{\lambda}^{-})}=0$, and
therefore $\hat{\Sigma}_{\lambda}^{-}$ must be measure zero. Hence, for
these values of $\lambda>\lambda_{0}$, we have
$$
w_{\lambda}(x)\geq 0,\;\;a.e.\;\forall x\in\hat{\Sigma}_{\lambda}\backslash B_{\epsilon}((z^{0})^{\lambda}),\;\;\forall\epsilon>0.
$$
This contradicts the definition of $\lambda_{0}$. Therefore
(\ref{b3}) must hold. That is, if $\lambda_{0}<z_{1}^{0} -\epsilon$, for any $\epsilon>0$, then we must have
$$
\bar{u}(x)\equiv \bar{u}_{\lambda_{0}}(x),\;\;a.e.\;\forall x\in \hat{\Sigma}_{\lambda_{0}}\backslash B_{\epsilon}((z^{0})^{\lambda_{0}}).
$$
Since $\bar{u}$ is singular at $z^{0}$, $\bar{u}$ must also be singular at $(z^{0})^{\lambda}$. This is impossible because $z^0$ is the only singularity of $\bar{u}$. Hence we must have $\lambda_0 \geq z_1^0 - \epsilon$. Since $\epsilon$ is
an arbitrary positive number, we have actually derived that
$$
w_{z_1^{0}}(x)\geq 0,\;\;a.e.\;\forall x\in\hat{\Sigma}_{z_1^{0}}.
$$

Entirely similarly,  we can move the plane from near $x_{1}=\infty$ to the left and derive that
$w_{z_1^{0}}(x)\leq 0$. Therefore we have
$$
w_{z_1^{0}}(x)\equiv 0,\;\;a.e.\;\forall x\in \hat{\Sigma}_{z_1^{0}}.
$$

Now we prove inequality (\ref{b30}). For any small $\eta>0$, $\forall \epsilon>0$, we
can choose $R$ sufficiently large so that
\begin{equation}\label{b97}
\left(\int_{(\mathbb{R}^n_+\backslash B_{\epsilon}(z^{0}))\backslash
B_{R}}\bar{u}^{\frac{2n}{n-\alpha}}(y)dy\right)^{\frac{\alpha}{n}}<\eta.
\end{equation}
We fix this $R$ and then show that the measure of
$\hat{\Sigma}_{\lambda}^{-}\cap B_{R}$ is sufficiently small for
$\lambda$ close to $\lambda_{0}$. By (\ref{14}), we have
\begin{equation}\label{b25}
w_{\lambda_{0}}(x)>0
\end{equation}
in the interior of $\hat{\Sigma}_{\lambda_{0}}\backslash B_{\epsilon}((z^{0})^{\lambda_{0}})$.

The rest is similar to the proof of (\ref{a30a}). We only need to use $\hat{\Sigma}_{\lambda}\backslash B_{\epsilon}((z^{0})^{\lambda})$ instead of $\Sigma_{\lambda}$ and $\hat{\Sigma}_{\lambda_{0}}\backslash B_{\epsilon}((z^{0})^{\lambda_{0}})$ instead of $\Sigma_{\lambda_{0}}$.

\textbf{(ii)} {\em The Subcritical Case} $1<p<\frac{n+\alpha}{n-\alpha}$. In this case, we only need to carry the method of moving planes on
$\bar{u} \equiv \bar{u}_{z^0}$ to show that it must be axially symmetric about the line passing through $z^0$ and parallel to $x_n$ axis.

Since
$u$ is locally bounded, for any domain $\Omega$
that is a positive distance away from $z^{0}$, we have
\begin{equation}\label{c85}
\int_{\Omega}[\frac{\bar{u}^{p-1}(y)}{|y-z^{0}|^{\beta}}]^{\frac{n}{\alpha}}dy<\infty.
\end{equation}
where $\beta=(n-\alpha)(\tau-p)>0$, $\tau=\frac{n+\alpha}{n-\alpha}$.

By (\ref{53}), we have
$$
\bar{u}(x)=\int_{\hat{\Sigma}_{\lambda}}G_{\infty}(x,y)\frac{\bar{u}^{p}(y)}{|y-z^{0}|^{\beta}}dy
+\int_{\hat{\Sigma}_{\lambda}}G_{\infty}(x,y^{\lambda})
\frac{\bar{u}^{p}_{\lambda}(y)}{|y^{\lambda}-z^{0}|^{\beta}}dy,
$$
$$
\bar{u}(x^{\lambda})
=\int_{\hat{\Sigma}_{\lambda}}G_{\infty}(x^{\lambda},y)\frac{\bar{u}^{p}(y)}{|y-z^{0}|^{\beta}}dy
+\int_{\hat{\Sigma}_{\lambda}}G_{\infty}(x^{\lambda},y^{\lambda})\frac{\bar{u}^{p}_{\lambda}(y)}{|y^{\lambda}-z^{0}|^{\beta}}dy.
$$
From (\ref{54}), we calculate
\begin{eqnarray}\nonumber
&&\bar{u}(x)-\bar{u}_{\lambda}(x)\\
\nonumber
&=&\int_{\hat{\Sigma}_{\lambda}}\left[G_{\infty}(x,y)-G_{\infty}(x^{\lambda},y)\right]
\frac{\bar{u}^{p}(y)}{|y-z^{0}|^{\beta}}dy  \\
\nonumber &+&\int_{\hat{\Sigma}_{\lambda}}
\left[G_{\infty}(x,y^{\lambda})-G_{\infty}
(x^{\lambda},y^{\lambda})\right]\frac{\bar{u}^{p}_{\lambda}(y)}{|y^{\lambda}-z^{0}|^{\beta}}dy\\
\label{c54}
&=&\int_{\hat{\Sigma}_{\lambda}}\left[
G_{\infty}(x,y)-G_{\infty}
(x^{\lambda},y)\right]\left[\frac{\bar{u}^{p}(y)}{|y-z^{0}|
^{\beta}}-\frac{\bar{u}_{\lambda}^{p}(y)}{|y^{\lambda}-z^{0}|^{\beta}}\right]dy.
\end{eqnarray}

The proof also consists of two steps.

\emph{Step 1.} For any $\epsilon>0$, define
$$
\hat{\Sigma}_{\lambda}^{-}=\{x\in \hat{\Sigma}_{\lambda}\backslash B_{\epsilon}((z^{0})^{\lambda})|\;\;w_{\lambda}(x)=\bar{u}_{\lambda}(x)-\bar{u}(x)<0\}.
$$
We show that for $\lambda$ sufficiently  negative,
$\hat{\Sigma}_{\lambda}^{-}$ must be measure zero.

By the Mean Value Theorem, we obtain, for sufficiently negative values of $\lambda$ and $x\in\hat{\Sigma}_{\lambda}^{-}$,
\begin{eqnarray}
\nonumber
0&<&\bar{u}(x)-\bar{u}_{\lambda}(x)\\
\nonumber
&= &\int_{\hat{\Sigma}_{\lambda}}\left[
G_{\infty}(x,y)-G_{\infty}
(x^{\lambda},y)\right]\left[\frac{\bar{u}^{p}(y)}{|y-z^{0}|^{\beta}}
-\frac{\bar{u}_{\lambda}^{p}(y)}{|y^{\lambda}-z^{0}|^{\beta}}\right]dy\\
\nonumber
&=&\int_{\hat{\Sigma}_{\lambda}^{-}}\left[
G_{\infty}(x,y)-G_{\infty}
(x^{\lambda},y)\right]\left[\frac{\bar{u}^{p}(y)}{|y-z^{0}|^{\beta}}
-\frac{\bar{u}_{\lambda}^{p}(y)}{|y^{\lambda}-z^{0}|^{\beta}}\right]dy\\
\nonumber
&\;\;+&\int_{\hat{\Sigma}_{\lambda}\backslash \hat{\Sigma}_{\lambda}^{-}}\left[
G_{\infty}(x,y)-G_{\infty}
(x^{\lambda},y)\right]\left[\frac{\bar{u}^{p}(y)}{|y-z^{0}|^{\beta}}
-\frac{\bar{u}_{\lambda}^{p}(y)}{|y^{\lambda}-z^{0}|^{\beta}}\right]dy\\
\nonumber
&\leq& \int_{\hat{\Sigma}_{\lambda}^{-}}\left[
G_{\infty}(x,y)-G_{\infty}
(x^{\lambda},y)\right]\left[\frac{\bar{u}^{p}(y)}{|y-z^{0}|^{\beta}}
-\frac{\bar{u}_{\lambda}^{p}(y)}{|y^{\lambda}-z^{0}|^{\beta}}\right]dy\\
\nonumber
&=&\int_{\hat{\Sigma}^{-}_{\lambda}}\left[
G_{\infty}(x,y)-G_{\infty}(x^{\lambda},y)\right]\left[\frac{\bar{u}^{p}(y)}{|y-z^{0}|^{\beta}}
-\frac{\bar{u}_{\lambda}^{p}(y)}{|y-z^{0}|^{\beta}}+\frac{\bar{u}_{\lambda}^{p}(y)}{|y-z^{0}|^{\beta}}-
\frac{\bar{u}_{\lambda}^{p}(y)}{|y^{\lambda}-z^{0}|^{\beta}}\right]dy\\
\nonumber
&=&\int_{\hat{\Sigma}^{-}_{\lambda}}\left[
G_{\infty}(x,y)-G_{\infty}(x^{\lambda},y)\right]\left[\frac{\bar{u}^{p}(y)-\bar{u}_{\lambda}^{p}(y)}
{|y-z^{0}|^{\beta}}+\bar{u}_{\lambda}^{p}(y)
[\frac{1}{|y-z^{0}|^{\beta}}-\frac{1}{|y^{\lambda}-z^{0}|^{\beta}}]\right]dy\\
\nonumber
&\leq&\int_{\hat{\Sigma}^{-}_{\lambda}}\left[
G_{\infty}(x,y)-G_{\infty}(x^{\lambda},y)\right]
\frac{\bar{u}^{p}(y)-\bar{u}_{\lambda}^{p}(y)}
{|y-z^{0}|^{\beta}}dy\\
\nonumber
&\leq&p\int_{\hat{\Sigma}^{-}_{\lambda}}
G_{\infty}(x,y)\frac{\bar{u}^{p-1}(y)}{|y-z^{0}|^{\beta}}[\bar{u}(y)-\bar{u}_{\lambda}(y)]dy\\
\label{c23}
&\leq&\int_{\hat{\Sigma}^{-}_{\lambda}}
\frac{C}{|x-y|^{n-\alpha}}|\frac{\bar{u}^{p-1}(y)}{|y-z^{0}|^{\beta}}||\bar{u}(y)-\bar{u}_{\lambda}(y)|dy.
\end{eqnarray}
We apply the Hardy-Littlewood-Sobolev inequality (\ref{f1}) and H\"{o}lder
inequality to (\ref{c23}) to obtain, for any $q>\frac{n}{n-\alpha}$,
\begin{eqnarray}\nonumber
\|w_{\lambda}\|_{L^{q}(\hat{\Sigma}_{\lambda}^{-})}&\leq&
C\|\frac{\bar{u}^{p-1}}{|y-z^{0}|^{\beta}}w_{\lambda}\|_{L^{\frac{nq}{n+\alpha
q}}(\hat{\Sigma}_{\lambda}^{-})}\\
\label{c24}
&\leq& C\|\frac{\bar{u}^{p-1}}{|y-z^{0}|^{\beta}}\|_{L^{\frac{n}{\alpha}}(\hat{\Sigma}_{\lambda}^{-})}
\|w_{\lambda}\|_{L^{q}(\hat{\Sigma}_{\lambda}^{-})}.
\end{eqnarray}
By (\ref{c85}), we can choose $N$ sufficiently large, such that for $\lambda\leq-N$,
\begin{equation}\label{h1}
C\left
\{\int_{\hat{\Sigma}_{\lambda}^{-}}
[\frac{\bar{u}^{p-1}}{|y-z^{0}|^{\beta}}]^{\frac{n}{\alpha}}dy\right\}^{\frac{\alpha}{n}}\leq\frac{1}{2}.
\end{equation}
Now inequality (\ref{c24}) and (\ref{h1}) imply
$$
\|w_{\lambda}\|_{L^{q}(\hat{\Sigma}_{\lambda}^{-})}=0,
$$
and therefore $\hat{\Sigma}_{\lambda}^{-}$ must be measure zero. Then we get
\begin{equation}\label{e11}
w_{\lambda}(x)\geq 0, \;\;a.e.\;x\in \hat{\Sigma}_{\lambda}.
\end{equation}

\emph{Step 2.} ({\em Move the plane to the limiting position to derive
symmetry}.)

Inequality (\ref{e11}) provides a starting point to move the plane
$\hat{T}_{\lambda}$. Now we start from
the neighborhood of $x_{1}=-\infty$ and move the plane to the right as
long as (\ref{e11}) holds to the limiting position. Define
$$
\lambda_{0}=\sup\{\lambda\leq z_{1}^{0}|\; w_{\rho}(x)\geq 0,\; \rho\leq\lambda,\;
\forall x\in\hat{\Sigma}_{\rho}\}.
$$
The rest is entirely similarly to the critical case when $p=\frac{n+\alpha}{n-\alpha}$. We only need to use $\int[\frac{\bar{u}^{p-1}(y)}{|y-z^{0}|^{\beta}}]^{\frac{n}{\alpha}}dy$ instead of $\int u^{\frac{2n}{n-\alpha}}(y)dy$. We also conclude
$$
w_{\lambda_{0}}(x)\equiv 0,\;\;a.e.\;\forall x\in \hat{\Sigma}_{\lambda_{0}},\;\;\lambda_{0}=z^{0}_{1}.
$$
This implies that $\bar{u}$ is symmetric about the plane $\hat{T}_{z^0}$.

Since we can choose any direction that is perpendicular to the $x_{n}$-axis as the $x_{1}$ direction, we have
actually shown that the Kelvin transform of the solution-- $\bar{u}(x)$-- is rotationally symmetric about the line parallel to $x_{n}$-axis and passing through $z^{0}$ either in \textbf{Possibility 2} of the critical case or in the subcritical case. Now, for any two points $X^{1}$ and $X^{2}$, with $X^{i}=(x^{i},x_{n})\in R^{n-1}\times [0,\infty)$, $i=1,2$.  Let $z^{0}$ be the projection of $\bar{X}=\frac{X^{1}+X^{2}}{2}$ on $\partial \mathbb{R}^n_+$. Set $Y^{i}=\frac{X^{i}-z_{0}}{|X^{i}-z^{0}|^{2}}+z^{0}$, $i=1,2$. From the above arguments, it is easy to see $\overline{u}(Y^{1})=\overline{u}(Y^{2})$, hence $u(X^{1})=u(X^{2}).$ This implies that $u$ is independent of $(x_{1},\cdots,x_{n-1})$. That is $u= u(x_n)$, and we will show that this will contradict the finiteness of the integral
$$ \int_{\mathbb{R}^n_+} G_{\infty} (x, y) u^p (y) d y .$$

Recall that
\begin{equation}\label{d1}
G_{\infty}(x,y)=\frac{A_{n,\alpha}}{s^{\frac{n-\alpha}{2}}}
\left[1-B\int_{0}^{\frac{s}{t}}\frac{
\left(\frac{1-\frac{t}{s}b}{1+\frac{t}{s}}\right)^{\frac{n-2}{2}}}{b^{\alpha/2}(1+b)}db\right]=\frac{A_{n,\alpha}}{s^{\frac{n-\alpha}{2}}}
\left[1-BI(s,t)\right],
\end{equation}
where $t=4x_{n}y_{n}$ and $s=|x-y|^{2}$. By the boundary conditions, we have
\begin{equation}
B=\frac{1}{\int_{0}^{\infty}\frac{1}{b^{\alpha/2}(1+b)}db}.
\end{equation}
For simplicity, we set $\lambda=\frac{t}{s}.$ Now we estimate
\begin{eqnarray}
1-BI(s,t)=\frac{\int_{0}^{\infty}\frac{1}{b^{\alpha/2}(1+b)}db-
\int_{0}^{\frac{1}{\lambda}}\frac{\left(\frac{1-\lambda b}{1+\lambda}\right)^{\frac{n-2}{2}}}{b^{\alpha/2}(1+b)}db}{\int_{0}^{\infty}\frac{1}{b^{\alpha/2}(1+b)}db}.
\end{eqnarray}
We calculate
\begin{eqnarray}
\nonumber
&&\int_{0}^{\infty}\frac{1}{b^{\alpha/2}(1+b)}db-
\int_{0}^{\frac{1}{\lambda}}\frac{\left(\frac{1-\lambda b}{1+\lambda}\right)^{\frac{n-2}{2}}}{b^{\alpha/2}(1+b)}db\\
&=&\int_{0}^{\frac{1}{\lambda}}\frac{1-\left(\frac{1-\lambda b}{1+\lambda}\right)^{\frac{n-2}{2}}}{b^{\alpha/2}(1+b)}db
+\int_{\frac{1}{\lambda}}^{\infty}\frac{1}{b^{\alpha/2}(1+b)}db.
\end{eqnarray}

Let $s\rightarrow \infty$, that is $\lambda\rightarrow 0.$  We have
\begin{equation}
\int_{\frac{1}{\lambda}}^{\infty}\frac{1}{b^{\alpha/2}(1+b)}db
\sim\int_{\frac{1}{\lambda}}^{\infty}\frac{1}{b^{\alpha/2+1}}db\sim\lambda^{\alpha/2}.
\end{equation}
Here we use ``$\sim$'' to indicate that the two quantities have the same order.

Let $b=\frac{1}{\lambda}\widetilde{b}$ and $s\rightarrow \infty$. Then we have
\begin{eqnarray}
\nonumber
\int_{0}^{\frac{1}{\lambda}}\frac{1-\left(\frac{1-\lambda b}{1+\lambda}\right)^{\frac{n-2}{2}}}{b^{\alpha/2}(1+b)}db&=&\left(\frac{1}{1+\lambda}\right)^{\frac{n-2}{2}}
\int_{0}^{1}\frac{\lambda^{\alpha/2}[(1+\lambda)^{\frac{n-2}{2}} -(1-\widetilde{b})
^{\frac{n-2}{2}}]}{{\widetilde{b}^{\alpha/2}(\lambda+\widetilde{b})}}d\widetilde{b}\\
\nonumber
&\sim&\frac{\lambda^{\alpha/2}}{(1+\lambda)^{\frac{n-2}{2}}}
\int_{0}^{1}\frac{(\lambda +\widetilde{b}) \eta^{\frac{n-2}{2}-1}
}{{\widetilde{b}^{\alpha/2}(\lambda+\widetilde{b})}}d\widetilde{b}\\
\label{d2}
&\sim& \lambda^{\alpha/2}.
\end{eqnarray}
Here we have used the Mean Value Theorem with $\eta$ valued between $1-\widetilde{b}$ and $1+\lambda$.
By (\ref{d1})-(\ref{d2}), we derive that, for each fixed $t>0$ and for sufficiently large $s$,
\begin{equation}
\frac{c_{n,\alpha}}{s^{\frac{n-\alpha}{2}}}\cdot\frac{t^{\alpha/2}}{s^{\alpha/2}}\leq G_{\infty}(x,y)\leq
\frac{C_{n,\alpha}}{s^{\frac{n-\alpha}{2}}}\cdot\frac{t^{\alpha/2}}{s^{\alpha/2}}.
\end{equation}
That is
\begin{equation}\label{11}
G_{\infty}(x,y)\sim\frac{t^{\alpha/2}}{s^{\frac{n}{2}}}.
\end{equation}
Set $x=(x',x_{n})$, $y=(y',y_{n})\in R^{n-1}\times(0,+\infty)$, $r^{2}=|x'-y'|^{2}$ and $a^{2}=|x_{n}-y_{n}|^{2}$.
If $u(x)=u(x_{n})$ is a solution of
\begin{equation}\label{w10}
u(x)=\int_{\mathbb{R}^n_+}G_{\infty}(x,y)u^{p}(y)dy,
\end{equation}
then for each fixed $x\in \mathbb{R}^n_+$, let $\rho, R>1$ be large enough, by (\ref{11}), we have
\begin{eqnarray}
\nonumber
+\infty>u(x_{n})&=&\int_{0}^{\infty}u^{p}(y_{n})\int_{R^{n-1}}G_{\infty}(x,y)dy'dy_{n}\\
\nonumber
&\geq&\int_{\rho}^{\infty}u^{p}(y_{n})\int_{R^{n-1}\backslash B_{R}(0)}G_{\infty}(x,y)dy'dy_{n}\\
\nonumber
&\geq &C\int_{\rho}^{\infty}u^{p}(y_{n})y_{n}^{\alpha/2}\int_{R^{n-1}\backslash B_{R}(0)}\frac{1}{|x-y|^{n}}dy'dy_{n}\\
\nonumber
&\geq &C\int_{\rho}^{\infty}u^{p}(y_{n})y_{n}^{\alpha/2}\int_{R}^{\infty}
\frac{r^{n-2}}{(r^{2}+a^{2})^{\frac{n}{2}}}drdy_{n}\\
\nonumber
&\geq&C\int_{\rho}^{\infty}u^{p}(y_{n})y_{n}^{\alpha/2}\frac{1}{|x_{n}-y_{n}|}\int_{R}^{\infty}
\frac{\tau^{n-2}}{(\tau^{2}+1)^{\frac{n}{2}}}d\tau dy_{n}\\
\label{12}
&\sim&C_{R}\int_{\rho}^{\infty}u^{p}(y_{n})y_{n}^{\alpha/2-1}dy_{n}.
\end{eqnarray}
(\ref{12}) implies that there exists a sequence $\{y_{n}^{i}\}\rightarrow\infty$ as $i\rightarrow\infty$, such that
\begin{equation}\label{w100}
u^{p}(y^{i}_{n})(y^{i}_{n})^{\alpha/2}\rightarrow 0.
\end{equation}

Similarly to (\ref{12}), for any $x=(0,x_{n})\in \mathbb{R}^n_+$, we derive that
\begin{equation}\label{w1}
+\infty>u(x_{n})\geq C_{0}\int_{0}^{\infty}u^{p}(y_{n})y_{n}^{\alpha/2}\frac{1}{|x_{n}-y_{n}|}dy_{n}x_{n}^{\alpha/2}.
\end{equation}
Let $x_{n}=2R$ be sufficiently large. By (\ref{w1}), we deduce that
\begin{eqnarray}
\nonumber
+\infty>u(x_{n})&\geq& C_{0}\int_{0}^{1}u^{p}(y_{n})y_{n}^{\alpha/2}\frac{1}{|x_{n}-y_{n}|}dy_{n}x_{n}^{\alpha/2} \\
\nonumber
&\geq& \frac{C_{0}}{2R}(2R)^{\alpha/2}\int_{0}^{1}u^{p}(y_{n})y_{n}^{\alpha/2}dy_{n}\\
\label{w2}
&\geq& C_{1}(2R)^{\alpha/2-1}=C_{1}x_{n}^{\alpha/2-1}.
\end{eqnarray}
Then by (\ref{w1}) and (\ref{w2}), for $x_{n}=2R$ sufficiently large, we also obtain
\begin{eqnarray}
\nonumber
u(x_{n})&\geq& C_{0}\int_{\frac{R}{2}}^{R}(C_{1}y_{n})^{p(\alpha/2-1)}y_{n}^{\alpha/2}\frac{1}{|x_{n}-y_{n}|}dy_{n}x_{n}^{\alpha/2} \\
\nonumber
&\geq& C_{0}(C_{1})^{p(\frac{\alpha}{2}-1)}R^{p(\frac{\alpha}{2}-1)}\frac{2}{3R}(2R)^{\alpha/2}\int_{\frac{R}{2}}^{R}y_{n}^{\alpha/2}dy_{n} \\
\nonumber
&\geq&  C_{0}(C_{1})^{p(\frac{\alpha}{2}-1)}\frac{2^{\alpha/2+2}}{3(\alpha+2)}(1-\frac{1}{2^{\alpha/2+1}})
R^{p(\frac{\alpha}{2}-1)+\alpha} \\
\nonumber
&:=&AR^{p(\frac{\alpha}{2}-1)+\alpha}\\
\nonumber
&=& \frac{A}{2^{p(\frac{\alpha}{2}-1)+\alpha}}x_{n}^{p(\frac{\alpha}{2}-1)+\alpha} \\
\label{w3}
&:=& A_{1}x_{n}^{p(\frac{\alpha}{2}-1)+\alpha}.
\end{eqnarray}
Continuing this way $m$ times, for $x_{n}=2R$, we have
\begin{equation}\label{w4}
u(x_{n})\geq A(m,p,\alpha)x_{n}^{p^{m}(\frac{\alpha}{2}-1)+\frac{p^{m}-1}{p-1}\alpha}.
\end{equation}
For any fixed $0<\alpha<2$, we choose $m$ to be an integer greater than $\frac{3-\alpha^{2}}{\alpha}.$ That is
\begin{equation}\label{w5}
m\geq \left\lceil\frac{3-\alpha^{2}}{\alpha}\right\rfloor+1,
\end{equation}
where $\lceil a\rfloor$ is the integer part of $a$.

We claim that for such choice of $m$, it holds
\begin{equation}\label{w6}
\tau (p):=\left[p^{m}(\frac{\alpha}{2}-1)+\frac{p^{m}-1}{p-1}\alpha\right]p+\frac{\alpha}{2}\geq 0.
\end{equation}
We postpone the proof of (\ref{w6}) for a moment. Now by (\ref{w4}) and (\ref{w6}), we derive that
\begin{equation}
u^{p}(x_{n})x_{n}^{\alpha/2}\geq  A(m,p,\alpha)x_{n}^{\tau (p)}\geq A(m,p,\alpha)>0,
\end{equation}
for all $x_{n}$ sufficiently large. This contradicts (\ref{w100}). So there is no positive solution of (\ref{w10}).

Now what left is to verify (\ref{w6}). In fact, if we let
$$
f(p):=\tau (p)(p-1)=p^{m+2}(\frac{\alpha}{2}-1)+(\frac{\alpha}{2}+1)p^{m+1}-\frac{\alpha}{2}p-\frac{\alpha}{2},
$$
then
$$
f'(p)=p^{m}[(m+2)(\frac{\alpha}{2}-1)p+(m+1)(\frac{\alpha}{2}+1)]-\frac{\alpha}{2}.
$$
We show that
$$
f'(p)>0,\;\;\mbox{for}\;\;1<p\leq\frac{n+\alpha}{n-\alpha}.
$$
Since $p>1$, it suffices to show
$$
(m+2)(\frac{\alpha}{2}-1)p+(m+1)(\frac{\alpha}{2}+1)\geq\frac{\alpha}{2}.
$$
Due to the fact $\frac{\alpha}{2}-1<0$, $n\geq 3,$ and $p\leq\frac{n+\alpha}{n-\alpha}$, we only
need to verify that
$$
(m+2)(\frac{\alpha}{2}-1)\frac{3+\alpha}{3-\alpha}+(m+1)(\frac{\alpha}{2}+1)\geq\frac{\alpha}{2}
$$
which can be derived directly from (\ref{w5}).

This completes the proof of Theorem \ref{mthm4}.

Authors' Addresses and Emails
\medskip

Wenxiong Chen

Department of Mathematics

Yeshiva University

2495 Amsterdam Av.

New York NY 10033, USA

wchen@yu.edu
\medskip

Yanqin Fang

School of Mathematics

Hunan University

Changsha, 410082 P.R. China, and

Department of Mathematics

Yeshiva University

yanqinfang@hnu.edu.cn
\medskip

Ray Yang

Department of Mathematics

Courant Institute of Mathematical Sciences

251 Mercer Street

New York NY 10012, USA

ryang@cims.nyu.edu

\end{document}